\input amstex
\input amsppt.sty \magnification=1200
\NoBlackBoxes
\hsize=15.6truecm \vsize=22.2truecm
\def\q{\quad}
\def\qq{\qquad}
\def\qtq#1{\q\t{#1}\q}

%\par\q
%\hbox{The paper is submitted to Advances in Applied Mathematics for
%possible publication.}\par\q \hbox{(January 18, 2012)}\par\q\par
\def\({\left(}
\def\){\right)}
\def\[{\left[}
\def\]{\right]}
\def\mod#1{\ (\text{\rm mod}\ #1)}
\def\t{\text}
\def\f{\frac}

\def\e{\equiv}
\def\a{\alpha}
\def\b{\binom}

\def\ap{\langle a\rangle_p}
\def\bp{\langle b\rangle_p}

\def\sls#1#2{(\f{#1}{#2})}
\def\ag#1{\langle #1\rangle_p}
\def\ls#1#2{\big(\f{#1}{#2}\big)}
\def\Ls#1#2{\Big(\f{#1}{#2}\Big)}
\let \pro=\proclaim
\let \endpro=\endproclaim

\topmatter
\title {Generalized Legendre polynomials and related congruences modulo $p^2$}\endtitle
\author Zhi-Hong Sun \endauthor
\affil School of Mathematical Sciences, Huaiyin Normal University,
\\ Huaian, Jiangsu 223001, PR China
\\ Email: zhihongsun$\@$yahoo.com
\\ Homepage: http://www.hytc.edu.cn/xsjl/szh
\endaffil

 \nologo \NoRunningHeads
%\TagsOnRight
 \nologo \NoRunningHeads

\abstract{For any positive integer $n$ and variables $a$ and $x$ we
define the generalized Legendre polynomial $P_n(a,x)=\sum_{k=0}^n\b
ak\b{-1-a}k(\frac{1-x}2)^k$. Let $p$ be an odd prime. In the paper
we prove many congruences modulo $p^2$ related to $P_{p-1}(a,x)$.
For example, we show that $P_{p-1}(a,x)\e (-1)^{\langle a\rangle
_p}P_{p-1}(a,-x)\mod {p^2}$, where $\langle a\rangle _p$ is the
least nonnegative residue of $a$ modulo $p$. We also generalize some
congruences of Zhi-Wei Sun, and determine
$\sum_{k=0}^{p-1}\binom{2k}k\binom{3k}k{54^{-k}}$ and
$\sum_{k=0}^{p-1}\binom ak\binom{b-a}k\mod {p^2}$, where $[x]$ is
the greatest integer function. Finally we pose some supercongruences
modulo $p^2$ concerning binary quadratic forms.
\par\q
}
 \endabstract
  \footnote"" {The author was
supported by the Natural Sciences Foundation of China (grant no.
10971078).}
\endtopmatter
\document
\subheading{1. Introduction}
\par
For any positive integer $n$ and variables $a$ and $x$ we introduce
the generalized Legendre polynomial
$$\aligned P_n(a,x)&=\sum_{k=0}^n\b
ak\b{-1-a}k\Ls{1-x}2^k=\sum_{k=0}^n\b ak\b{a+k}k\Ls{x-1}2^k
\\&=\sum_{k=0}^n\b{a+k}{2k}\b{2k}k\Ls{x-1}2^k.\endaligned\tag 1.1$$
 We note that
$\b{-1-a}k=(-1)^k\b{a+k}k$ and $\b ak\b{a+k}k=\b{a+k}{2k}\b{2k}k$.
Clearly $P_n(a,x)=P_n(-1-a,x)$ and $P_n(n,x)=P_n(x)$ (see [B]),
where $P_n(x)$ is the Legendre polynomial given by
$$P_n(x)=\f
1{2^n}\sum_{k=0}^{[n/2]}\b nk(-1)^k\b{2n-2k}nx^{n-2k} =\f 1{2^n\cdot
n!}\cdot\f{d^n}{dx^n}(x^2-1)^n.\tag 1.2$$
\par Let $p>3$ be a prime. In 2003,
Rodriguez-Villegas[RV] conjectured the following congruences:
$$\align &\sum_{k=0}^{p-1}\f{\b{2k}k^2}{16^k}\e
\Ls{-1}p\mod{p^2},\tag 1.3
\\&\sum_{k=0}^{p-1}
\f{\b{2k}k\b{3k}k}{27^k}\e\Ls {-3}p\mod{p^2},\tag 1.4
\\&\sum_{k=0}^{p-1}\f{\b{2k}k\b{4k}{2k}}{64^k}\e \Ls{-2}p\mod{p^2},\tag 1.5
\\& \sum_{k=0}^{p-1} \f{\b{3k}k\b{6k}{3k}}{432^k}\e
\Ls{-1}p\mod{p^2},\tag 1.6\endalign$$ where $\sls ap$ is the
Legendre symbol. These congruences were later confirmed by
Mortenson[M1-M3] via the Gross-Koblitz formula.
 Recently Zhi-Wei Sun[Su1] posed more conjectures concerning the
following sums modulo $p^2$:
$$\sum_{k=0}^{p-1}\f{\b{2k}k^2}{16^k}x^k,\ \sum_{k=0}^{p-1}
\f{\b{2k}k\b{3k}k}{27^k}x^k,\
\sum_{k=0}^{p-1}\f{\b{2k}k\b{4k}{2k}}{64^k}x^k,\ \sum_{k=0}^{p-1}
\f{\b{3k}k\b{6k}{3k}}{432^k}x^k.$$ As observed by Tauraso[T] and
Zudilin, we have
$$\align&\b{-\f 12}k^2=\f{\b{2k}k^2}{16^k},\ \b{-\f
13}k\b{-\f 23}k=\f{\b{2k}k\b{3k}k}{27^k}, \\&\ \b{-\f 14}k\b{-\f
34}k=\f{\b{2k}k\b{4k}{2k}}{64^k},\ \b{-\f 16}k\b{-\f
56}k=\f{\b{3k}k\b{6k}{3k}}{432^k}.\endalign$$
 This is the motivation that we
introduce and study $P_{p-1}(a,x)\mod {p^2}$.
\par Let $\Bbb Z$ be the set of integers. For a prime $p$ let
$\Bbb Z_p$ denote the set of rational $p$-integers. For a $p$-integr
$a$ let $\ap$ be the least nonnegative residue of $a$ modulo $p$.
Let $p$ be an odd prime $p$ and $a\in\Bbb Z_p$. In the paper we show
that
$$P_{p-1}(a,x)\e (-1)^{\ap}P_{p-1}(a,-x)\mod{p^2}.\tag 1.7$$
Taking $x=-1$ we obtain
$$\sum_{k=0}^{p-1}\b ak\b{-1-a}k=P_{p-1}(a,-1)\e (-1)^{\ap}\mod{p^2}.
\tag 1.8$$ For $a=-\f 12,-\f 13,-\f 14,-\f 16$ we get (1.3)-(1.6)
immediately from (1.8). If $\ap$ is odd, by (1.7) we have
$$\sum_{k=0}^{p-1}\b ak\b{-1-a}k\f 1{2^k}=
P_{p-1}(a,0)\e 0\mod{p^2}.\tag 1.9$$ This generalizes previous
special results in [S2] and [Su2]. If $f(0),f(1),\ldots,f(p-1)$ are
$p$-integers, we also prove the following general congruence:
$$\sum_{k=0}^{p-1}\b ak\b{-1-a}k\Big((-1)^{\ap}f(k)-\sum_{m=0}^k\b
km(-1)^mf(m)\Big)\e 0\mod{p^2}.\tag 1.10$$ (1.7)-(1.9) can be viewed
as vast generalizations of some congruences proved in [S2,S3] (with
$a=-\f 12$) and [Su2,Su3] (with $a=-\f 13,-\f 14,-\f 16$).
 When
$\langle a\rangle_p\e 1\mod 2$, taking $f(k)=\b{2k}k/2^{2k}$ in
(1.10) we get
$$\sum_{k=0}^{p-1}\b ak\b{-1-a}k\b{2k}k\f 1{4^k}\e 0\mod{p^2}.\tag 1.11$$
Inspired by the work of Zhi-Wei Sun([Su2-Su4]) we establish the
general congruence
$$(a+1)P_n(a+1,x)-(2a+1)xP_n(a,x)+aP_n(a-1,x)\e
0\mod{p^2}\ \t{for}\ a\not\e 0,-1\mod p.$$ and use it to prove our
main results. As an application, we deduce the congruence for
$\sum_{k=0}^{p-1}\b{2k}k\b{3k}k{54^{-k}}\mod{p^2}$, see Theorem 3.4.
In Section 4, we establish a general congruence for
$\sum_{k=0}^{p-1}\b ak\b{b-a}k\mod{p^2}$ and pose some conjectures
for $\sum_{k=0}^{p-1}\b ak\b{b-a}km^k\mod{p^2}$, where $p$ is an odd
prime and $a,b\in\Bbb Z_p$.

\subheading{2. General congruences for $P_{p-1}(a,x)\mod{p^2}$}
 \pro{Lemma 2.1} Let $n$ be a positive
integer. Then
$$\align&(a+1)P_n(a+1,x)-(2a+1)xP_n(a,x)+aP_n(a-1,x)
\\&=-2(2a+1)\b an\b{a+n}n\Ls{x-1}2^{n+1}.\endalign$$
\endpro
Proof. It is clear that
$$a(a+1-k)(a-k)+(a+1)(a+1+k)(a+k)
=(2a+1)((a+k)(a-k+1)+2k^2).$$ Thus,
$$\align &a\b{a-1+k}{2k}+(a+1)\b{a+1+k}{2k}-(2a+1)\b{a+k}{2k}-(2a+1)\f
k{2k-1}\b{a+k-1}{2k-2}
\\&=\f{(a+k-1)(a+k-2)\cdots (a-(k-2))}{(2k)!}
\big(a(a+1-k)(a-k)\\&\qq+(a+1)(a+1+k)(a+k)-(2a+1)((a+k)(a-k+1)-2k^2)\big)
\\&=0.\endalign$$
Therefore,
$$\aligned&a\b{a-1+k}{2k}\b{2k}k+(a+1)\b{a+1+k}{2k}\b{2k}k-(2a+1)\b{a+k}{2k}\b{2k}k
\\&=(2a+1)\cdot 2\b{a+k-1}{2k-2}\b{2k-2}{k-1}.\endaligned\tag 2.1$$
Using (2.1) we deduce that
$$\align&(a+1)P_n(a+1,x)-(2a+1)xP_n(a,x)+aP_n(a-1,x)
\\&=\sum_{k=0}^n\Big\{(a+1)\b{a+1+k}{2k}\b{2k}k+a\b{a-1+k}{2k}\b{2k}k
\\&\qq-(2a+1)\Big(1+2\cdot
\f{x-1}2\Big)\b{a+k}{2k}\b{2k}k\Big\}\Ls{x-1}2^k
\\&=-(2a+1)\cdot 2\b{a+n}{2n}\b{2n}n\Ls{x-1}2^{n+1}
\\&\q+\sum_{k=0}^n\Big\{(a+1)\b{a+1+k}{2k}\b{2k}k+a\b{a-1+k}{2k}\b{2k}k
\\&\qq-(2a+1)\Big(\b{a+k}{2k}\b{2k}k+2\b{a+k-1}{2k-2}\b{2k-2}{k-1}\Big)\Big\}
\Ls{x-1}2^k
\\&=-2(2a+1)\b an\b{a+n}n\Ls{x-1}2^{n+1}.
\endalign$$
This proves the lemma.

 \pro{Theorem 2.1} Let $p$ be an odd prime and $a\in\Bbb Z_p$.
Then
$$\aligned &(a+1)P_{p-1}(a+1,x)-(2a+1)xP_{p-1}(a,x)+aP_{p-1}(a-1,x)
\\&\e \cases -2p^2(\f 1{\ap}+\f 1{\ap+1})\f{a-\ap}p(1+\f{a-\ap}p)\sls{x-1}2^p
\mod{p^3}&\t{if $a\not\e 0,-1\mod p$,}
\\-2a(p+a+1)\sls{x-1}2^p \mod{p^3}&\t{if $a\e 0\mod p$,}
\\2(a+1)(a-p)\sls{x-1}2^p\mod{p^3}&\t{if $a\e -1\mod p$.}\endcases\endaligned$$
\endpro
Proof. Clearly
$$\b a{p-1}\b{a+p-1}{p-1}=
\f{a(a-1)\cdots (a-(p-2))(a+1)(a+2)\cdots (a+p-1)}{(p-1)!^2}.$$
 If
$\ap\not=0,p-1$, then
$$\align &a(a-1)\cdots
(a-(p-2))(a+1)(a+2)\cdots (a+p-1)
\\&\e (a-\ap)(a+p-\ap)\{\ap(\ap-1)\cdots 2\cdot 1\cdot
(p-1)(p-2)\cdots (\ap+2)
\\&\qq\times (\ap+1)(\ap+2)\cdots(p-1)\cdot 1\cdot
2\cdots(\ap-1)\}
\\&=p^2\cdot\f{a-\ap}p\Big(1+\f{a-\ap}p\Big)\f{(p-1)!^2}{\ap(\ap+1)}
\mod{p^3}.\endalign$$ If $\ap=0$, then
$$\align \b a{p-1}\b{a+p-1}{p-1}
&=a(a+p-1)\f{(a^2-1^2)(a^2-2^2)\cdots (a^2-(p-2)^2)}{(p-1)!^2}
\\&\e a(a+p-1)\f{(-1^2)(-2^2)\cdots (-(p-2)^2)}{(p-1)!^2}
\e -\f{a(a+p-1)}{(p-1)^2}
\\&\e -a(a+p-1)(p+1)^2\e -a(a+p-1)(2p+1)
\\&\e -a(-2p+a+p-1)=a(p+1-a)\mod{p^3}.
\endalign$$
If $\ap=p-1$, then $-1-a\e 0\mod p$ and so
$$\align\b a{p-1}\b{a+p-1}{p-1}&=\b{-1-a}{p-1}\b{-1-a+p-1}{p-1}
\\&\e (-1-a)(p+1-(-1-a))=-(a+1)(p+a+2)\mod{p^3}.\endalign$$ Now putting
all the above together with Lemma 2.1 in the case $n=p-1$ we deduce
the result. \pro{Corollary 2.1} Let $p$ be an odd prime and
$a\in\Bbb Z_p$ with $a\not\e 0,-1\mod p$. Let $P_n'(a,x)=\f
d{dx}P_n(a,x)$. Then
$$\align&(a+1)P'_{p-1}(a+1,x)-(2a+1)xP'_{p-1}(a,x)-(2a+1)P_{p-1}(a,x)+aP'_{p-1}(a-1,x)
\\&\e 0\mod{p^2}.\endalign$$
\endpro
Proof. By Theorem 2.1 we have
$$(a+1)P_{p-1}(a+1,x)-(2a+1)xP_{p-1}(a,x)+aP_{p-1}(a-1,x)=p^2f(x),$$
where $f(x)$ is a polynomial of degree at most $p-1$ with rational
$p$-integral coefficients. Taking derivatives on both sides we
deduce the result.

 \pro{Lemma 2.2} Let
$p$ be an odd prime and let $t$ and $x$ be $p$-adic integers. Then
$$P_{p-1}(pt,x)\e 1-t+t\Ls{1+x}2^p+t\Ls{1-x}2^p\mod{p^2}\qtq{and so}
P_{p-1}(pt,x)\e 1\mod p.$$
\endpro
Proof. It is clear that
$$\align P_{p-1}(pt,x)&=\sum_{k=0}^{p-1}\b{pt}k\b{pt+k}k\Ls{x-1}2^k
\\&=1+t\sum_{k=1}^{p-1}\f pk\b{pt-1}{k-1}\b{pt+k}k\Ls{x-1}2^k
\\&\e 1+t\sum_{k=1}^{p-1}\b
pk\Ls{x-1}2^k=1+t\Big\{\Big(1+\f{x-1}2\Big)^p-1-\Ls{x-1}2^p\Big\}
\\&=1-t+t\Big\{\Ls{1+x}2^p+\Ls{1-x}2^p\Big\}\mod{p^2}.\endalign$$
Since $\ls{1\pm x}2^p\e \f{1\pm x^p}2\mod p$, by the above we obtain
$P_{p-1}(pt,x)\e 1\mod p$. This proves the lemma.

 \pro{Lemma 2.3} Let $p$ be an odd prime and let $t$ and $x$ be $p$-adic integers. Then
$$\align &P_{p-1}(1+pt,x)
\\&\e(1-t)x+pt\f{x-1}2+t\Big\{(p+1)\f{x-1}2\Big(\Ls{1+x}2^p+\Ls{1-x}2^p\Big)
\\&\qq+(p+1)^2\Ls{1-x}2^p+\Ls{1+x}2^{p+1}-\Ls{1-x}2^{p+1}\Big\}\mod{p^2}
\endalign$$ and so
$P_{p-1}(1+pt,x)\e x\mod p$.
\endpro
Proof. It is clear that
$$\align &P_{p-1}(1+pt,x)\\&=\sum_{k=0}^{p-1}\b{pt+1}k\b{pt+1+k}k\Ls{x-1}2^k
\\&=1+(pt+1)(pt+2)\f{x-1}2+\sum_{k=2}^{p-1}\f {pt+1}k\cdot \f{pt}{k-1}\b{pt-1}{k-2}\b{pt+1+k}k\Ls{x-1}2^k
\\&\e 1+(3pt+2)\f{x-1}2+t\sum_{k=2}^{p-1}\b
{p+1}k\b{k+1}k\Ls{x-1}2^k
\\&=1+(pt+2-2t)\f{x-1}2+t\sum_{k=1}^{p-1}\b{p+1}k(k+1)\Ls{x-1}2^k\mod{p^2}
\endalign$$ and so
$$\align&P_{p-1}(1+pt,x)
\\&\e 1+(pt+2-2t)\f{x-1}2\\&\qq+t\sum_{k=1}^{p-1}\b{p+1}k\Ls{x-1}2^k
+t(p+1)\f{x-1}2\sum_{k=1}^{p-1}
\b p{k-1}\Ls{x-1}2^{k-1}
\\&=1+(pt+2-2t)\f{x-1}2+t\Big\{\Big(1+\f{x-1}2\Big)^{p+1}-1-(p+1)\Ls{x-1}2^p
-\Ls{x-1}2^{p+1}
\\&\qq+(p+1)\f{x-1}2\Big(\Big(1+\f{x-1}2\Big)^p-\Ls{x-1}2^p-p\Ls{x-1}2^{p-1}\Big)\Big\}
\\&=(1-t)x+pt\f{x-1}2+t\Big\{(p+1)\f{x-1}2\Big(\Ls{1+x}2^p+\Ls{1-x}2^p\Big)
\\&\qq+(p+1)^2\Ls{1-x}2^p+\Ls{1+x}2^{p+1}-\Ls{1-x}2^{p+1}\Big\}\mod{p^2}.
\endalign$$
Observe that $\ls{1\pm x}2^p\e \f{1\pm x^p}2\mod p$. From the above
we see that
$$\align P_{p-1}(1+pt,x)&\e
(1-t)x+t\Big\{\f{x-1}2+\f{1-x^p}2+\f{1+x}2\cdot\f{1+x^p}2-\f{1-x}2\cdot\f{1-x^p}2\Big\}
\\&=x-tx+tx=x\mod p.\endalign$$ This completes the proof.

\pro{Theorem 2.2} Let $p$ be an odd prime and $a\in\Bbb Z_p$. Then
$$P_{p-1}(a,x)\e (-1)^{\langle a\rangle_p}P_{p-1}(a,-x)\mod{p^2}$$
and so
$$\sum_{k=0}^{p-1}\b ak\b{-1-a}k(x^k-(-1)^{\langle a\rangle_p}(1-x)^k)\e
0\mod{p^2}.$$
\endpro
Proof. Suppose $m\in\{1,2,\ldots,p-2\}$ and $t\in\Bbb Z_p$. From
Theorem 2.1 we have
$$\align (m+1&+pt)P_{p-1}(m+1+pt,\pm x)-(2(m+pt)+1)(\pm x)P_{p-1}(m+pt,\pm x)
\\&\qq+(m+pt)P_{p-1}(m-1+pt,\pm x)\e 0\mod{p^2}.\endalign$$
Thus,
$$\aligned &\q(m+1+pt)(P_{p-1}(m+1+pt,x)-(-1)^{m+1}P_{p-1}(m+1+pt,-x))
\\&\qq\e (2(m+pt)+1)x(P_{p-1}(m+pt,x)-(-1)^mP_{p-1}(m+pt,-x))
\\&\qq\q-(m+pt)(P_{p-1}(m-1+pt,x)-(-1)^{m-1}P_{p-1}(m-1+pt,-x))
\mod{p^2}.\endaligned\tag 2.2$$ From Lemma 2.2 we know that
$$P_{p-1}(pt,x)\e P_{p-1}(pt,-x)\mod{p^2}.$$
From Lemma 2.3 we see that
$$\align &P_{p-1}(1+pt,x)+P_{p-1}(1+pt,-x)
\\&\e pt\Big(\f{x-1}2+\f{-x-1}2\Big)+t\Big\{(p+1)\Big(\f{x-1}2+\f{-x-1}2\Big)
\Big(\Ls{1+x}2^p+\ls{1-x}2^p\Big)
\\&\qq+(1+2p)\Big(\Ls{1+x}2^p+\ls{1-x}2^p\Big)\Big\}
\\&\e -pt+pt\Big(\Ls{1+x}2^p+\ls{1-x}2^p\Big)\e
-pt+pt\Big(\f{1+x^p}2 +\f{1-x^p}2\Big)
\\&=0\mod{p^2}.\endalign$$
Thus,
$$P_{p-1}(m+pt,x)-(-1)^mP_{p-1}(m+pt,-x)\e
0\mod{p^2}\qtq{for}m=0,1.$$ By (2.2) and induction we deduce that
$P_{p-1}(m+pt,x)-(-1)^mP_{p-1}(m+pt,-x)\e 0\mod{p^2}$ for all
$m=0,1,\ldots,p-1$. Thus, $P_{p-1}(a,x)\e (-1)^{\langle
a\rangle_p}P_{p-1}(a,-x)$ and hence
$$\align \sum_{k=0}^{p-1}\b ak\b{-1-a}kx^k&=P_{p-1}(a,1-2x)\e(-1)^{\langle
a\rangle_p}P_{p-1}(a,2x-1)\\&=(-1)^{\langle
a\rangle_p}\sum_{k=0}^{p-1}\b
ak\b{-1-a}k(1-x)^k\mod{p^2}.\endalign$$ This completes the proof.
\par\q
\newline{\bf Remark 2.1} In the case $a=-\f 12$, Theorem 2.2 was given by
the author in [S2]. In the cases $a=-\f 13,-\f 14,-\f 16$, Theorem
2.2 was given by Z. W. Sun in [Su3].

\pro{Corollary 2.2} Let $p$ be an odd prime and $a\in\Bbb Z_p$. Then
$$P_{p-1}(a,-1)=\sum_{k=0}^{p-1}\b ak\b{-1-a}k\e (-1)^{\ap}\mod{p^2}.$$
\endpro
Proof. Taking $x=1$ in Theorem 2.2 we obtain the result.
\par\q As mentioned in Section 1, taking $a=-\f 12,-\f 13,-\f 14,-\f 16$ in Corollary 2.2
we deduce (1.3)-(1.6).

\pro{Corollary 2.3} Let $p$ be an odd prime and $a\in\Bbb Z_p$ with
$\langle a\rangle_p\e 1\mod 2$. Then
$$\sum_{k=0}^{p-1}\b ak\b{-1-a}k\f 1{2^k}\e 0\mod{p^2}.$$
\endpro
Proof. Taking $x=\f 12$ in Theorem 2.2 we obtain the result.
\par\q Putting $a=-\f 12,-\f 13,-\f 14,-\f 16$ in Corollary 2.3
we deduce the following congruences
$$\align &\sum_{k=0}^{p-1}\f{\b{2k}k^2}{32^k}\e
0\mod{p^2}\qtq{for}p\e 3\mod 4,\tag 2.3
\\&\sum_{k=0}^{p-1}
\f{\b{2k}k\b{3k}k}{54^k}\e 0\mod{p^2}\qtq{for}p\e 2\mod 3,\tag 2.4
\\&\sum_{k=0}^{p-1}\f{\b{2k}k\b{4k}{2k}}{128^k}\e
0\mod{p^2}\qtq{for}p\e 5,7\mod 8,\tag 2.5
\\& \sum_{k=0}^{p-1} \f{\b{3k}k\b{6k}{3k}}{864^k}\e
0\mod{p^2}\qtq{for}p\e 3\mod 4,\tag 2.6\endalign$$ where $p$ is a
prime greater than $3$. We remark that (2.3) was conjectured by Z.
W. Sun and proved by the author in [S2], and (2.4) was conjectured
by the author in [S2] and proved by Z. W. Sun in [Su4]. (2.5) and
(2.6) were conjectured by Z. W. Sun and finally proved by him in
[Su4], although the author proved the congruences modulo $p$
earlier.
 \pro{Lemma
2.4} Let $n$ be a positive integer. Then
$$\sum_{k=0}^n\f{\b ak\b{-1-a}k}{k+1}=\f{\b{a-1}n\b{-2-a}n}{n+1}.$$
\endpro
Proof. Set
$$f(n)=\f{\b
an\b{-1-a}n}{n+1}\qtq{and}g(n)=\f{\b{a-1}n\b{-2-a}n}{n+1}.$$ It is
easily seen that $g(n)-g(n-1)=f(n)$. Thus,
$$\sum_{k=0}^nf(k)=f(0)+\sum_{k=1}^n(g(k)-g(k-1))
=f(0)-g(0)+g(n)=g(n).$$ This proves the lemma.

 \pro{Theorem 2.3} Let $p$ be an odd prime, $a\in\Bbb
Z_p$ and let $m$ be a positive integer with $m<p$. Then
$$\sum_{k=m}^{p-1}\b ak\b{-1-a}k\b km(x^{k-m}-(-1)^{m+\ap}(1-x)^{k-m})\e
0\mod{p^2}$$ and
$$\aligned&\sum_{k=0}^{p-2}\b ak\b{-1-a}k\f{x^{k+1}+(-1)^{\ap}(1-x)^{k+1}}{k+1}
\\&\e\cases 0\mod{p^2}&\t{if $a\not\e 0,-1\mod p$,}
\\1-2a-a(x^p-(x-1)^p-1)/p\mod{p^2}&\t{if $a\e 0\mod p$,}
\\2a+3+(1+a)(x^p-(x-1)^p-1)/p\mod{p^2}&\t{if $a\e -1\mod p$.}\endcases
\endaligned$$
\endpro
Proof. By Theorem 2.2,
$$\sum_{k=0}^{p-1}\b ak\b{-1-a}k(x^k-(-1)^{\langle a\rangle_p}(1-x)^k)=p^2f(x),$$
where $f(x)$ is a polynomial of $x$ with rational $p$-integral
coefficients and degree at most $p-1$. Since $\f
{d^m\;x^k}{dx^m}=m!\b kmx^{k-m}$ and
$\f{d^m}{dx^m}(1-x)^k=(-1)^km!\b km(x-1)^{k-m}=(-1)^mm!\b
km(1-x)^{k-m}$ for $k\ge m$, we see that
$$\sum_{k=m}^{p-1}\b ak\b{-1-a}k\b km(x^{k-m}-(-1)^{m+\ap}(1-x)^{k-m})
=\f {p^2}{m!}\cdot \f{d^m}{dx^m}f(x).$$ As $\f{d^m}{dx^m}f(x)$ is a
polynomial of $x$ with rational $p$-integral coefficients and degree
at most $p-1$, we deduce the first result.
\par Set
$$\delta_p(a)=\cases 0&\t{if $a\not\e 0,-1\mod
p$,}\\a&\t{if $a\e 0\mod p$,}
\\-1-a&\t{if $a\e -1\mod p$.}
\endcases$$ It is easy to see that
$$\b a{p-1}\b{-1-a}{p-1}\e \delta_p(a)\mod{p^2}.$$
Thus, by Theorem 2.2 we have
$$\sum_{k=0}^{p-2}\b ak\b{-1-a}k(x^k-(-1)^{\ap}(1-x)^k)
+\delta_p(a)(x^{p-1}-(1-x)^{p-1})=p^2g(x),$$ where $g(x)$ is a
polynomial of $x$ with rational $p$-integral coefficients and degree
at most $p-1$. It is easy to see that
$$\align &p^2\int_0^xg(t)dt\\&=\sum_{k=0}^{p-2}\b ak\b{-1-a}k\f{x^{k+1}+(-1)^{\ap}((1-x)^{k+1}-1)}
{k+1}+\delta_p(a)\f{x^p-(x-1)^p-1}p.\endalign$$ Thus, using  Lemma
2.4 we get $$\align&\sum_{k=0}^{p-2}\b
ak\b{-1-a}k\f{x^{k+1}+(-1)^{\ap}(1-x)^{k+1}}{k+1}
\\&\e \sum_{k=0}^{p-2}\f{\b ak\b{-1-a}k}{k+1}
-\delta_p(a)\f{x^p-(x-1)^p-1}p
\\&=\f 1{p-1}\b{a-1}{p-2}\b{-2-a}{p-2}-\delta_p(a)\f{x^p-(x-1)^p-1}p
\mod{p^2}.\endalign$$ If $a\not\e 0,-1\mod p$, then $\delta_p(a)=0$
and
$$\b{a-1}{p-2}\b{-2-a}{p-2}=\f a{p-1}\cdot \f{p-1}{-1-a}\b a{p-1}\b{-1-a}{p-1}
\e 0\mod{p^2}.$$ If $a\e 0\mod p$, then $\delta_p(a)=a$ and
$$\align\f{\b{a-1}{p-2}\b{-2-a}{p-2}}{p-1}&=\f{(a-1)(a-2)\cdots
(a-(p-2))(-2-a)\cdots(-(p-1)-a)}{(p-2)!(p-1)!}\\& \e
\f{(p-2)!(1-a\sum_{i=1}^{p-2}\f 1i)\cdot
(p-1)!(1+a\sum_{i=2}^{p-1}\f 1i)}{(p-2)!(p-1)!}
\\&\e (1+a/(p-1))(1-a)\e 1-2a\mod{p^2}.
\endalign$$
If $a\e -1\mod p$, then $\delta_p(a)=-1-a$. Using the above we get
$$\align\f 1{p-1}\b{a-1}{p-2}\b{-2-a}{p-2}&=\f
1{p-1}\b{-1-a-1}{p-2}\b{-2-(-1-a)}{p-2}\\&\e
1-2(-1-a)=2a+3\mod{p^2}.\endalign$$ Now putting all the above
together we deduce the result.
 \pro{Corollary 2.4} Let $p$ be an odd
prime, $a\in\Bbb Z_p$ and let $m$ be a positive integer with $m<p$.
Then
$$\sum_{k=m}^{p-1}\b ak\b{-1-a}k\b km\e (-1)^{m+\langle a\rangle_p}\b
am\b{-1-a}m\mod{p^2}.$$
\endpro
Proof. Taking $x=1$ in Theorem 2.3 we obtain the result.

 \pro{Theorem 2.4} Let $p$ be an odd
prime and $a\in\Bbb Z_p$. If $f(0),f(1),\ldots,f(p-1)$ are
$p$-integers, then
$$\sum_{k=0}^{p-1}\b ak\b{-1-a}k\Big((-1)^{\ap}f(k)-\sum_{m=0}^k\b
km(-1)^mf(m)\Big)\e 0\mod{p^2}.$$
\endpro
Proof. Using Corollary 2.4 we see that
$$\align&\sum_{k=0}^{p-1}\b ak\b{-1-a}k\sum_{m=0}^k\b km(-1)^mf(m)
\\&=\sum_{m=0}^{p-1}(-1)^mf(m)\sum_{k=m}^{p-1}\b ak\b{-1-a}k\b km
\\&\e \sum_{m=0}^{p-1}(-1)^mf(m)\cdot (-1)^{m+\ap}\b am\b{-1-a}m
\\&=(-1)^{\ap}\sum_{k=0}^{p-1}\b ak\b{-1-a}kf(k)\mod{p^2}.\endalign$$
This yields the result.
\newline{\bf Remark 2.2} In the case $a=-\f 12$,
 Corollary 2.3 and Theorem 2.4 were obtained by the author in
[S3]. In the cases $a=-\f 13,-\f 14,-\f 16$, Corollary 2.3 and
Theorem 2.4 were recently obtained by Z.W. Sun in [Su3].
\pro{Theorem 2.5} Let $p$ be an odd prime and $a\in\Bbb Z_p$ with
$\ap\e 1\mod 2$. Then
$$\sum_{k=0}^{p-1}\b ak\b{-1-a}k\b{2k}k\f 1{2^{2k}}\e 0\mod
{p^2}.$$
\endpro
Proof. Set $f(k)=\b{2k}k/2^{2k}$. From [S1, Example 10] we know that
$\sum_{m=0}^k\b km(-1)^mf(m)=f(k)$. Thus, applying Theorem 2.4 we
deduce the result.
\par Putting $a=-\f 12,-\f 13,-\f 14,-\f 16$ in Theorem 2.5 we deduce that
for any prime $p>3$,
$$\align &\sum_{k=0}^{p-1}\f{\b {2k}k^3}{64^k}\e 0\mod{p^2}\qtq{for}
p\e 3\mod 4,\tag 2.7
\\&\sum_{k=0}^{p-1}\f{\b{2k}k^2\b{3k}k}{108^k}\e 0\mod{p^2}\qtq{for}
p\e 5\mod 6,\tag 2.8
\\&\sum_{k=0}^{p-1}\f{\b{2k}k^2\b{4k}{2k}}{256^k}\e 0\mod{p^2}\qtq{for}
p\e 5,7\mod 8,\tag 2.9
\\&\sum_{k=0}^{p-1}\f{\b{2k}k\b{3k}k\b{6k}{3k}}{432^k}\e 0\mod{p^2}\qtq{for}
p\e 3\mod 4.\tag 2.10\endalign$$ Here (2.7) was conjectured by
Beukers[Be] in 1985 and proved by van Hamme[vH]. (2.8)-(2.10) were
conjectured by Rodriguez-Villegas[RV] and proved by Z. W. Sun[Su4].

\subheading{3. Congruences for $P_{p-1}(a,0)\mod{p^2}$}
\par For given positive integer $n$ and prime $p$ we define
$$H_n=1+\f 12+\f 13+\cdots+\f 1n\qtq{and}q_p(a)=\f{a^{p-1}-1}p.$$
 \pro{Theorem 3.1} Let $p$ be an odd prime and $t\in\Bbb Z_p$.
\par $(\t{\rm i})$ If $n\in\{0,1,\ldots,\f{p-3}2\}$, then
$P_{p-1}(2n+1+pt,0)\e 0\mod{p^2}.$
\par $(\t{\rm ii})$ If $n\in\{0,1,\ldots,\f{p-1}2\}$, then
 $$P_{p-1}(2n+pt,0)\e \b{\f{p-1}2}n\Big(1+p((1+t)H_{2n}-\f
{2t+1}2H_n-tq_p(2))\Big)\mod{p^2}.$$
\endpro
Proof. Putting $x=0$ in Theorem 2.1 we see that
$$P_{p-1}(a+1,0)\e -\f a{a+1}P_{p-1}(a-1,0)\mod{p^2}\qtq{for}a\not\e
0,-1\mod p.$$ Assume $n\in\{1,2,\ldots,\f{p-3}2\}$. Then
$$\align &P_{p-1}(2n+1+pt,0)\\&=-\f{2n+pt}{2n+1+pt}P_{2n-1+pt}(0)=\cdots
\\&=(-1)^n\f{2n+pt}{2n+1+pt}\cdot
\f{2n-2+pt}{2n-1+pt}\cdots\f{2+pt}{3+pt}P_{p-1}(1+pt,0) \mod{p^2}.
\endalign$$
By Lemma 2.3 we have
$$P_{p-1}(1+pt,0)\e -\f p2t+t\Big(-\f{p+1}2\cdot \f
1{2^{p-1}}+(2p+1)\cdot \f 1{2^p}\Big)=-\f p2t+\f p{2^p}t\e
0\mod{p^2}.$$ Thus, from the above we deduce that
$P_{p-1}(2n+1+pt,0)\e 0\mod{p^2}$ for $n=0,1,\ldots,\f{p-3}2$. This
proves (i).
\par Now let us consider (ii). Assume $n\in\{1,2,\ldots,\f{p-1}2\}$.
Then
$$\align &P_{p-1}(2n+pt,0)\\&=-\f{2n-1+pt}{2n+pt}P_{2n-2+pt},0)=\cdots
\\&=(-1)^n\f{2n-1+pt}{2n+pt}\cdot
\f{2n-3+pt}{2n-2+pt}\cdots\f{3+pt}{4+pt}\cdot\f{1+pt}{2+pt}P_{p-1}(pt,0)
\\&\e (-1)^n\f{1\cdot 3\cdots (2n-1)(1+pt\sum_{k=1}^n\f 1{2k-1})}{2\cdot 4\cdots (2n)
(1+pt\sum_{k=1}^n\f 1{2k})} P_{p-1}(pt,0)
\\&\e\f 1{(-4)^n}\b{2n}n\Big(1+pt\sum_{k=1}^n\f 1{2k-1}\Big)
\Big(1-pt\sum_{k=1}^n\f 1{2k}\Big)P_{p-1}(pt,0)
\\&\e \f 1{(-4)^n}\b{2n}n(1+pt(H_{2n}-H_n))P_{p-1}(pt,0)
\mod{p^2}.\endalign$$ By [S1, Lemma 2.4] we have
$$\f 1{(-4)^n}\b{2n}n\e \b{\f{p-1}2}n\Big(1+p\sum_{k=1}^n\f 1{2k-1}\Big)
=\b{\f{p-1}2}n\Big(1+p(H_{2n}-\f 12H_n)\Big)\mod{p^2}.\tag 3.1$$
From Lemma 2.2 we have
$$P_{p-1}(pt,0)\e 1-t+\f t{2^{p-1}}\e
1-t+t(1-pq_p(2))=1-ptq_p(2)\mod{p^2}.$$ Therefore, for
$n=0,1,\ldots,\f{p-1}2$ we have
$$\align P_{p-1}(2n+pt,0)&\e \b{\f{p-1}2}n(1+p(H_{2n}-\f
{H_n}2))(1+pt(H_{2n}-H_n))(1-ptq_p(2))
\\&\e\b{\f{p-1}2}n(1+p((1+t)H_{2n}-\f{2t+1}2H_n))(1-ptq_p(2))
\\&\e \b{\f{p-1}2}n(1+p((1+t)H_{2n}-\f{2t+1}2H_n-tq_p(2)))\mod{p^2}.
\endalign$$ This proves (ii) and hence the proof is complete.
\pro{Corollary 3.1} Let $p$ be an odd prime and let $a$ be a
$p$-adic integer with $a\not\e 0\mod p$. Then $P_{p-1}(a,0)\e
0\mod{p^2}$ or $P_{p-1}(-a,0)\e 0\mod{p^2}$.
\endpro
Proof. If $a\e 2n\mod p$ for some $n\in\{1,2,\ldots,\f{p-1}2\}$,
then $-a\e p-2n=1+2(\f{p-1}2-n)\mod p$. Thus the result follows from
Theorem 3.1(i).

 \pro{Theorem 3.2} Let $p$ be an odd prime,
$r,m\in\Bbb Z$, $m\ge 1$, $r\in\{\pm 1,\pm 2,\ldots,\pm m-1\}$ and
$(r,m)=1$.
\par $(\t{\rm i})$ If $2\nmid rm$, then
$$\aligned P_{p-1}\Big(\f rm,0\Big)&=\sum_{k=0}^{p-1}\b {r/m}k\b{-1-r/m}k\f 1{2^k}
\\&\e\cases \b{(p-1)/2}n(1+p((1-\f{2s}m)H_{2n}+(\f{2s}m-\f 12)H_n+\f
{2s}mq_p(2)))\mod {p^2} \\\qq\t{if $p\e\f r{2s}\mod m$ for some
$s\in\{1,2,\ldots,\f{m-1}2\}$}\\\qq\t{with $(s,m)=1$ and
$n=\f{sp-(m+r)/2}m$,}
\\0\mod{p^2}\qq\t{otherwise.}
\endcases\endaligned$$
\par $(\t{\rm ii})$ If $2\nmid m$ and $2\mid r$, then
$$\aligned P_{p-1}\Big(\f rm,0\Big)&=\sum_{k=0}^{p-1}\b {r/m}k\b{-1-r/m}k\f 1{2^k}
\\&\e\cases  \b{(p-1)/2}n(1+p((1-\f{2s}m)H_{2n}+(\f{2s}m-\f 12)H_n+\f
{2s}mq_p(2)))\mod {p^2}\\\qq\t{if $p\e-\f {r/2}s\mod m$ for some
$s\in\{1,2,\ldots,\f{m-1}2\}$}\\\qq\t{with $(s,m)=1$ and
$n=\f{sp+r/2}m$,}
\\0\mod{p^2}\qq\t{otherwise.}
\endcases\endaligned$$
\par $(\t{\rm iii})$ If $2\mid m$ and $2\nmid r$, then
$$\aligned P_{p-1}\Big(\f rm,0\Big)&=\sum_{k=0}^{p-1}\b {r/m}k\b{-1-r/m}k\f 1{2^k}
\\&\e\cases \b{(p-1)/2}n(1+p((1-\f{s}m)H_{2n}+(\f{s}m-\f 12)H_n+\f
{s}mq_p(2)))\mod {p^2} \\\qq\t{if $p\e-\f rs\mod {2m}$ for some
$s\in\{1,2,\ldots,\f m2\}$}\\\qq\t{with $(s,m)=1$ and
$n=\f{sp+r}{2m}$,}
\\0\mod{p^2}\qq\t{otherwise.}
\endcases\endaligned$$
\endpro
Proof. If $2\nmid rm$, setting $n=\f{sp-(m+r)/2}m$ and $t=-\f{2s}m$
we find $2n+pt=-1-\f rm$. If $2\nmid m$ and $2\mid r$, setting
$n=\f{sp+r/2}m$ and $t=-\f{2s}m$ we find $2n+pt=\f rm$. If $2\mid m$
and $2\nmid r$, setting $n=\f{sp+r}{2m}$ and $t=-\f sm$ we find
$2n+pt=\f rm$. Now applying Theorem 3.1 and the fact $P_{p-1}(\f
rm,x)=P_{p-1}(-1-\f rm,x)$ we deduce the result.
\par From Theorem 3.2 we deduce the following result.
\pro{Theorem 3.3} Let $p$ be an odd prime. Then
$$\align &P_{p-1}(1/2,0)\e 0\mod{p^2}\qtq{for}p\e 1\mod 4,
\\&P_{p-1}(1/3,0)\e 0\mod{p^2}\qtq{for}p\e 1\mod 3,
\\&P_{p-1}(1/4,0)\e 0\mod{p^2}\qtq{for}p\e 1,3\mod 8,
\\&P_{p-1}(1/5,0)\e 0\mod{p^2}\qtq{for}p\e 1,2\mod 5,
\\&P_{p-1}(1/6,0)\e 0\mod{p^2}\qtq{for}p\e 1\mod 4,
\\&P_{p-1}(1/7,0)\e 0\mod{p^2}\qtq{for}p\e 1,3,5\mod 7,
\\&P_{p-1}(1/8,0)\e 0\mod{p^2}\qtq{for}p\e 1,7,11,13\mod {16},
\\&P_{p-1}(1/9,0)\e 0\mod{p^2}\qtq{for}p\e 1,2,4\mod 9,
\\&P_{p-1}(1/10,0)\e 0\mod{p^2}\qtq{for}p\e 1,3,7,9\mod {20},
\\&P_{p-1}(1/11,0)\e 0\mod{p^2}\qtq{for}p\e 1,4,5,8,9\mod {11},
\\&P_{p-1}(1/12,0)\e 0\mod{p^2}\qtq{for}p\e 1,5,7,11\mod {24}.
\endalign$$
\endpro

\pro{Lemma 3.1} Let $p$ be an odd prime. Then
\par $(\t{\rm i})\ (\t{\rm [L]})$ $H_{\f {p-1}2}\e -2q_p(2)\mod p$, $H_{[\f p4]}\e -3q_p(2)\mod p$.
\par $(\t{\rm ii})\ (\t{\rm [L]})$ For $p>3$ we have
$H_{[\f p3]}\e -\f 32q_p(3)\mod p$ and $H_{[\f p6]}\e -2q_p(2)-\f
32q_p(3)\mod p.$\endpro

\pro{Theorem 3.4} Let $p>3$ be a prime. Then
$$\sum_{k=0}^{p-1}\f{\b{2k}k\b{3k}k}{54^k}
\e\cases 2A-\f p{2A}\mod {p^2} &\t{if $3\mid p-1$, $p=A^2+3B^2$ and
$3\mid A-1$,}
\\ 0\mod {p^2}&\t{if $3\mid p-2$.}
\endcases$$\endpro
Proof. From Theorem 3.2(i) we see that
$$\aligned&\sum_{k=0}^{p-1}\f{\b{2k}k\b{3k}k}{54^k}=P_{p-1}\Big(-\f 13,0\Big)
\\&\e\cases \b{\f{p-1}2}{\f{p-1}3}(1+p(\f 13H_{\f{2(p-1)}3}+\f
16H_{\f{p-1}3}+\f 23q_p(2)))\mod{p^2}&\t{if $p\e 1\mod 3$,}
\\0\mod{p^2}&\t{if $p\e 2\mod 3$.}
\endcases\endaligned$$
\par Now we assume $p\e 1\mod 3$ and $p=A^2+3B^2$ with $A,B\in\Bbb
Z$ and $A\e 1\mod 3$. From Lemma 3.1 we have $H_{\f{p-1}3}\e -\f
32q_p(3)\mod p$ and
$$H_{\f{2(p-1)}3}=H_{p-1}-\sum_{k=1}^{(p-1)/3}\f 1{p-k}
\e \sum_{k=1}^{(p-1)/3}\f 1k\e -\f 32q_p(3)\mod p.$$ Thus,
$$\align\f 13H_{\f{2(p-1)}3}+\f
16H_{\f{p-1}3}+\f 23q_p(2)&\e -\f 12q_p(3)-\f 14q_p(3)+\f
23q_p(2)\\&=\f 23q_p(2)-\f 34q_p(3)\mod p.\endalign$$ By [BEW,
Theorem 9.4.4],
$$\b{\f{p-1}2}{\f{p-1}3}=\b{\f{p-1}2}{\f{p-1}6}
\e\Big(2A-\f p{2A}\Big)\Big(1-p\f 23q_p(2)+p\f
34q_p(3)\Big)\mod{p^2}.$$ Therefore,
$$\align P_{p-1}\Big(-\f 13,0\Big)&\e
\b{\f{p-1}2}{\f{p-1}3}\Big(1+p\big(\f 13H_{\f{2(p-1)}3}+\f
16H_{\f{p-1}3}+\f 23q_p(2)\big)\Big)
\\&\e\Big(2A-\f p{2A}\Big)\Big(1-p\Big(\f 23q_p(2)-\f
34q_p(3)\Big)\Big)\Big(1+p\Big(\f 23q_p(2)-\f 34q_p(3)\Big)\Big)
\\&\e 2A-\f p{2A}\mod {p^2}.\endalign$$
This completes the proof.
\par\q\newline{\bf Remark 3.1} In [S1] the author conjectured Theorem 3.4
and proved the congruence modulo $p$. In [Su3], Zhi-Wei Sun proved
the result for $p\e 2\mod 3$.

\subheading{4. Congruences for $\sum_{k=0}^{p-1}\b
ak\b{b-a}k\mod{p^2}$}
\par Let $n$ be a nonnegative integer. For two variables $a$ and $b$ we
define
$$S_n(a,b)=\sum_{k=0}^n\b ak\b{b-a}k.\tag 4.1$$
Using Maple and Zeilberger's algorithm we find
$$(a-b)S_n(a,b)+(a+1)S_n(a+1,b)=(2a-b+1)\b an\b{b-a-1}n.\tag 4.2$$

\pro{Lemma 4.1} Let $p$ be an odd prime and $b,c,t\in\Bbb Z_p$. Then
$$\sum_{k=0}^{p-1}\b{pt}k\b{b+cpt}k\e 1+pt\sum_{k=1}^{\langle
b\rangle_p}\f 1k\mod{p^2}.$$
\endpro
Proof. Clearly
$$\align &\sum_{k=0}^{p-1}\b{pt}k\b{b+cpt}k
\\&=1+\sum_{k=1}^{p-1}\f {pt}k\b{pt-1}{k-1}\b{b+cpt}k
\e 1+\sum_{k=1}^{p-1}\f{pt}k\b{-1}{k-1}\b bk
\\&\e 1+\sum_{k=1}^{p-1}\f {pt}k(-1)^{k-1}\b {\bp}k
=1+pt\sum_{k=1}^{\bp}\f{(-1)^{k-1}}k\b{\bp}k
\\&=1+pt\sum_{k=1}^{\bp}(-1)^{k-1}\b{\bp}k\int_0^1x^{k-1}\;dx
=1+pt\int_0^1\sum_{k=1}^{\bp}\b{\bp}k(-x)^{k-1}\;dx
\\&=1+pt\int_0^1\f{(1-x)^{\bp}-1}{-x}\;
dx=1+pt\int_0^1\f{u^{\bp}-1}{u-1}du
\\&=1+pt\int_0^1\sum_{k=0}^{\bp-1}u^k\;du
=1+pt\sum_{k=0}^{\bp-1}\f 1{k+1}\mod{p^2}.
\endalign$$
This proves the lemma.

\pro{Lemma 4.2} Let $p$ be an odd prime, $m\in\{1,2,\ldots,p-1\}$
and $t\in\Bbb Z_p$. Then
$$\b{m+pt-1}{p-1}\e
 \f{pt}m-\f{p^2t^2}{m^2}+\f{p^2t}m\sum_{k=1}^m\f 1k
\mod{p^3}.$$
\endpro
Proof. For $m<\f p2$ we see that
$$\align &\b{m+pt-1}{p-1}\\&=\f{(m-1+pt)(m-2+pt)\cdots (1+pt)\cdot
pt(pt-1)\cdots(pt-(p-1-m))}{(p-1)!}
\\&=\f{pt(p^2t^2-1^2)\cdots (p^2t^2-(m-1)^2)(pt-m)\cdots
(pt-(p-1-m))}{(p-1)!}
\\&\e pt \f{(m-1)!(-1)(-2)\cdots \cdots (-(p-1-m))}{(p-1)!}
\Big(1-pt\sum_{k=m}^{p-1-m}\f 1k\Big)
\\&=pt\cdot \f{(-1)^{p-1-m}\cdot (m-1)!}{(p-m)\cdots (p-1)}
\Big(1-pt\big(H_{p-1-m}-H_m+\f 1m\big)\Big) \mod{p^3}.\endalign$$
 For $m>\f p2$ we also have
$$\align &\b{m+pt-1}{p-1}\\&=\f{(m-1+pt)(m-2+pt)\cdots (1+pt)\cdot
pt(pt-1)\cdots(pt-(p-1-m))}{(p-1)!}
\\&=\f{pt(p^2t^2-1^2)\cdots (p^2t^2-(p-1-m)^2)(pt+p-m)\cdots
(pt+m-1)}{(p-1)!}
\\&\e pt \f{(-1)^{p-1-m}(p-1-m)!(m-1)!}{(p-1)!}
\Big(1+pt\sum_{k=p-m}^{m-1}\f 1k\Big)
\\&=pt\cdot \f{(-1)^{p-1-m}\cdot (m-1)!}{(p-m)\cdots (p-1)}
\Big(1-pt\big(H_{p-1-m}-H_m+\f 1m\big)\Big)\mod{p^3}.
\endalign$$
Since $(p-m)\cdots (p-1)\e (-1)^mm!(1-pH_m)\mod{p^2}$, by the above
we get
 $$\align \b{m+pt-1}{p-1}&\e pt\cdot \f{(-1)^m\cdot
(m-1)!}{(-1)^m\cdot m!(1-pH_m)}
\Big(1-\f{pt}m-pt(H_{p-1-m}-H_m)\Big)
\\&\e \f{pt}m(1+pH_m)\Big(1-\f{pt}m-pt(H_{p-1-m}-H_m)\Big)
\\&\e \f{pt}m-\f{p^2t}m\Big(\f tm+t(H_{p-1-m}-H_m)-H_m\Big)
 \mod{p^3}.\endalign$$
To see the result, we note that
$$\align H_{p-1-m}-H_m&=\sum_{r=1}^{p-1-m}\f 1r-\sum_{k=1}^m\f 1k
=\sum_{k=m+1}^{p-1}\f 1{p-k}-\sum_{k=1}^m\f 1k
\\&\e -\sum_{k=m+1}^{p-1}\f 1k-\sum_{k=1}^m\f 1k
=-\sum_{r=1}^{(p-1)/2}\Big(\f 1r+\f 1{p-r}\Big)
\\&\e 0\mod p.\endalign$$

 \pro{Lemma 4.3} Let $p$ be an odd prime, $m\in\{1,2,\ldots,p-1\}$ and
$b,t\in\Bbb Z_p$. Then
$$\b{b-(m+pt)}{p-1}\e\cases 1\mod{p^2}&\t{if $m\e b+1\mod p$,}
\\\f{b-\bp-pt}{m-1-b}\mod{p^2}&\t{if $m\le \bp\mod p$,}
\\\f{b-\bp-p(t+1)}{m-1-b}\mod{p^2}&\t{if $m>\bp+1\mod p$.}
\endcases$$
\endpro
Proof. If $m\e b+1\mod p$, setting $b-m=rp+p-1$ we find $r\in\Bbb
Z_p$ and so
$$\align\b{b-m-pt}{p-1}&=
\f{(p-1+p(r-t))(p-2+p(r-t))\cdots (1+p(r-t))}{(p-1)!}\\&\e
1+p(r-t)\sum_{k=1}^{p-1}\f 1k=1+p(r-t)\sum_{k=1}^{\f{p-1}2}\big(\f
1k+\f 1{p-k}\big)\e 1\mod{p^2}.\endalign$$ Now we assume $m\not\e
b+1\mod p$ and $b-m=\langle b-m\rangle_p+pr$. Then $r\in\Bbb Z_p$
and
$$\align &\b{b-(m+pt)}{p-1}\\&=\f{(\langle b-m\rangle_p+p(r-t))(\langle b-m\rangle_p-1+p(r-t))
\cdots (\langle
b-m\rangle_p-p+2+p(r-t))}{(p-1)!}
\\&\e \langle b-m\rangle_p!\cdot p(r-t)\cdot (p-1)(p-2)\cdots (2+\langle
b-m\rangle_p)
\\&=\f{p(r-t)\cdot (p-1)!}{1+\langle b-m\rangle_p}
\e \f{b-m-\langle b-m\rangle_p-pt}{m-1-b} \mod{p^2}.\endalign$$
Since $$ \langle b-m\rangle_p=\cases \bp-m&\t{if $m\le \bp$,}
\\p+\bp-m&\t{if $m>\bp+1$,}
\endcases$$
we see that
$$b-m- \langle b-m\rangle_p=\cases b-m-(\bp-m)=b-\bp&\t{if $m\le
\bp$,}\\b-m-(p+\bp-m)=b-\bp-p&\t{if $m>\bp+1$.}
\endcases$$
Now combining all the above we obtain the result.

 \pro{Theorem 4.1} Let $p$ be an odd prime and
$a,b\in\Bbb Z_p$. Then
$$\aligned\sum_{k=0}^{p-1}\b ak\b{b-a}k
\e\cases \f{(-1)^{\ap-\bp-1}}{(\ap-\bp)\b{\ap}{\bp}}(b-\bp) \mod
{p^2}\q\q\ \q\t{if $\ap>\bp$,}
\\\b{\bp}{\ap}(1+(b-\bp)H_{\bp}-(a-\ap)H_{\ap}\\\q-(b-a-\langle
b-a\rangle_p)H_{\langle b-a\rangle_p})\mod{p^2}\q\t{if $\ap\le
\bp$.}
\endcases\endaligned$$
\endpro
Proof. If $\ap=0$, the result follows from Lemma 4.1. From now on we
assume $\ap\ge 1$. Set $a=\ap+pt$ and $b=\bp+ps$. Then $s,t\in\Bbb
Z_p$. For $m\in\{1,2,\ldots,p-1\}$, by (4.2) and Lemmas 4.2-4.3 we
obtain
$$\aligned&(m+pt)S_{p-1}(m+pt,b)+(m+pt-b-1)S_{p-1}(m+pt-1,b)
\\&=(2(m+pt)-1-b) \b{m+pt-1}{p-1}\b{b-(m+pt)}{p-1}\\&\e\cases
(2m-1-b)\f{pt}m+p^2t(\f tm+H_m)\mod{p^3}&\t{if $m=\bp+1$,}
\\(2m-1-b)\f{pt}m\cdot \f{b-\bp-pt}{m-1-b}\mod{p^3}&\t{if $m\le
\bp$,}
\\(2m-1-b)\f{pt}m\cdot \f{b-\bp-p(t+1)}{m-1-b}\mod{p^3}&\t{if
$m>\bp+1$.}\endcases
\endaligned\tag 4.3$$ Hence, if $1\le
\ap\le \bp$, then
$$\align &S_{p-1}(a,b)\\&=S_{p-1}(\ap+pt,b)
\e -\f{\ap+pt-b-1}{\ap+pt}S_{p-1}(\ap-1+pt,b)
\\&=\f{1+\bp-\ap+p(s-t)}{\ap+pt}S_{p-1}(\ap-1+pt,b)
\\&\e \f{1+\bp-\ap+p(s-t)}{\ap+pt}\cdot \f{1+\bp-(\ap-1)+p(s-t)}{\ap-1+pt}
S_{p-1}(\ap-2+pt,b)
\\&\e \cdots\e
\prod_{k=1}^{\ap}\f{1+\bp-k+p(s-t)}{k+pt}\cdot
S_{p-1}(pt,b)\mod{p^2}.\endalign$$ Now applying Lemma 4.1 we see
that for $1\le \ap\le \bp$,
$$\align &S_{p-1}(a,b)\\&\e \f{\bp(\bp-1)\cdots(\bp-\ap+1)(1+p(s-t)
(H_{\bp}-H_{\langle
b-a\rangle_p}))}{\ap!(1+ptH_{\ap})}\big(1+ptH_{\bp}\big)
\\&\e \b{\bp}{\ap}(1+p(s-t)
(H_{\bp}-H_{\langle b-a\rangle_p}))(1-ptH_{\ap})(1+ptH_{\bp})
\\&\e \b{\bp}{\ap}(1+psH_{\bp}-ptH_{\ap}-(ps-pt)H_{\langle
b-a\rangle_p})
\\&=\b{\bp}{\ap}(1+(b-\bp)H_{\bp}-(a-\ap)H_{\ap}-(b-a-\langle
b-a\rangle_p)H_{\langle b-a\rangle_p}) \mod{p^2}.
\endalign$$
\par Now we assume $\ap>\bp$. Clearly
$$S_{p-1}(\bp+pt,b)=1+\sum_{k=1}^{p-1}\b{\bp+pt}k\f{b-\bp-pt}k
\b{b-\bp-pt-1}{k-1}\e 1\mod p.$$ By (4.2) and (4.3) we have
$$\align &(\bp+1+pt)S_{p-1}(\bp+1+pt,b)\\&\e -(\bp+pt-b)S_{p-1}(\bp+pt,b)
+(2m-1-b)\f{pt}m
\\&\e b-\bp-pt+pt=b-\bp\mod{p^2}.\endalign$$
Thus,
$$S_{p-1}(\bp+1+pt,b)\e \f{b-\bp}{\bp+1+pt}\e \f{b-\bp}{\bp+1}\mod{p^2}.
\tag 4.4$$ Now from (4.2)-(4.4) we deduce that for $\ap>\bp+1$,
$$\align S_{p-1}(a,b)&=S_{p-1}(\ap+pt,b)\e
\Big(\f{b+1}{\ap+pt}-1\Big)S_{p-1}(\ap+pt-1,b)
\\&\e
\Big(\f{b+1}{\ap+pt}-1\Big)\Big(\f{b+1}{\ap-1+pt}-1\Big)S_{p-1}(\ap-2+pt,b)
\\&\e \cdots\e \prod_{k=\bp+2}^{\ap}\Big(\f{b+1}{k+pt}-1\Big)\cdot
S_{p-1}(\bp+1+pt,b)
\\&\e \prod_{k=\bp+2}^{\ap}\f{\bp+1-k}{k}\cdot
\f{b-\bp}{\bp+1}
\\&=\f{(-1)^{\ap-\bp-1}}{(\ap-\bp)\b{\ap}{\bp}}(b-\bp)
\mod{p^2}.\endalign$$ This completes the proof.

\pro{Corollary 4.1} Let $p$ be an odd prime and $a,b\in\Bbb Z_p$.
Then
$$\sum_{k=0}^{p-1}\b ak\b{b-a}k\e \b{\bp}{\ap}\mod p.$$
\endpro

  \pro{Corollary 4.2}
Let $p$ be an odd prime and  $a\in\Bbb Z_p$. Then
$$\sum_{k=0}^{p-1}\b ak\b{-a}k\e\cases 0\mod{p^2}&\t{if
$a\not\e 0\mod p$,}
\\1\mod {p^2}&\t{if $a\e 0\mod p$.}
\endcases$$
\endpro
Proof. Taking $b=0$ in Theorem 4.1 we deduce the result.

\pro{Corollary 4.3} Let $p$ be an odd prime and $a\in\Bbb Z_p$. Then
$$\sum_{k=0}^{p-1}\b ak\b{1-a}k\e\cases 0\mod{p^2}&\t{if
$a(1-a)\not\e 0\mod p$,}
\\1+a\mod {p^2}&\t{if $a\e 0\mod p$,}
\\2-a\mod {p^2}&\t{if $a\e 1\mod p$.}
\endcases$$
\endpro
 Proof. Taking $b=1$ in Theorem 4.1 we deduce the result.
 \par We note that
 it can be easily proved that
 $$\align&\sum_{k=0}^n\b ak\b{-a}k=\b{n+a}n\b{n-a}n,\tag 4.5
 \\&\sum_{k=0}^n\b ak\b{1-a}k=-\f {a^2-a-n}{n^2}\b{a-2}{n-1}\b{-a-1}{n-1}
 .\tag 4.6\endalign$$
\pro{Theorem 4.2} Let $p$ be an odd prime. Then
$$\aligned&\sum_{k=0}^{p-1}\b{-\f 14}k\b{-\f 12}k
\\&\e \cases (-1)^{\f{p+1}4}p\b{\f{(p-1)}2}{\f{p+1}4}^{-1}\mod{p^2}&\t{if $p\e
3\mod 4$,}
\\(-1)^{\f{p-1}4}(2x-\f p{2x})\mod{p^2}&\t{if $p=x^2+y^2\e 1\mod 4$
and $4\mid x-1$.}
\endcases\endaligned$$
\endpro
Proof. Set $a=-\f 14$ and $b=-\f 34$. Then
$$\ap=\cases \f{p-1}4&\t{if $p\e 1\mod 4$,}
\\\f{3p-1}4&\t{if $p\e 3\mod 4$,}
\endcases\qtq{and}
\bp=\cases \f{3p-3}4&\t{if $p\e 1\mod 4$,}
\\\f{p-3}4&\t{if $p\e 3\mod 4$.}
\endcases$$
If $p\e 3\mod 4$, then $\ap>\bp$. Thus, by Theorem 4.1 we have
$$\aligned\sum_{k=0}^{p-1}\b{-\f 14}k\b{-\f 12}k
\e \f{(-1)^{\f{3p-1}4-\f{p-3}4-1}}
{(\f{3p-1}4-\f{p-3}4)\b{\f{3p-1}4}{\f{p-3}4}} \Big(-\f
34-\f{p-3}4\Big)
\e\f{3p}2\b{\f{3(p+1)}4}{\f{p+1}4}^{-1}\mod{p^2}.\endaligned$$ Since
$$\align \b{\f{3(p+1)}4}{\f{p+1}4}&=\f{(p-\f{p-3}4)(p-\f{p+1}4)\cdots
(p-\f{p-3}2)}{\f{p+1}4!}
\\&\e (-1)^{\f{p+1}4}\f{\f{p-3}2\cdots\f{p-3}4}{\f{p+1}4!}
=(-1)^{\f{p+1}4}\b{\f{p-1}2}{\f{p+1}4}\cdot\f{(p-3)/4}{(p-1)/2} \\&
\e \f 32(-1)^{\f{p+1}4}\b{\f{p-1}2}{\f{p+1}4}\mod p,
\endalign$$ by the above we obtain the result in the case $p\e 3\mod 4$.

\par Now we assume $p\e 1\mod 4$ and so $p=x^2+y^2$ with $x,y\in\Bbb
Z$ and $x\e 1\mod 4$. By the proof of Lemma 4.2 and Lemma 3.1(i) we
have $H_{\f{3(p-1)}4}\e H_{\f{p-1}4}\e -3q_p(2)\mod p$ and
$H_{\f{p-1}2}\e -2q_p(2)$. Now applying the above and Theorem 4.1 we
deduce
$$\aligned\sum_{k=0}^{p-1}\b{-\f 14}k\b{-\f 12}k
&\e\b{\f{3(p-1)}4}{\f{p-1}4}\Big(1-\f{3p}4H_{\f{3(p-1)}4}+\f
p4H_{\f{p-1}4}+\f p2H_{\f{p-1}2}\Big)
\\&\e \b{\f{3(p-1)}4}{\f{p-1}4}\Big(1-\f p2(-3q_p(2))+\f
p2(-2q_p(2))\Big)
\\&=\b{\f{3(p-1)}4}{\f{p-1}4}\Big(1+\f 12pq_p(2)\Big)\mod{p^2}.
\endaligned$$
By [BEW, Theorem 9.4.3] we have
$$\b{\f{3(p-1)}4}{\f{p-1}4}\e \Big(2x-\f
p{2x}\Big)(-1)^{\f{p-1}4}\Big(1-\f 12pq_p(2)\Big)\mod{p^2}.$$ Hence
$$\sum_{k=0}^{p-1}\b{-\f 14}k\b{-\f 12}k\e (-1)^{\f{p-1}4}\Big(2x-\f
p{2x}\Big)\mod{p^2}.$$ This proves the result in the case $p\e 1\mod
4$. The proof is now complete.

\pro{Theorem 4.3} Let $p>3$ be a prime.  Then
$$\aligned&\sum_{k=0}^{p-1}\b{-\f 16}k\b{-\f 13}k
\\&\e \cases \f{3p}2\b{\f{p-1}2}{\f{p-5}6}^{-1}\mod
{p^2}&\t{if $p\e 2\mod 3$,}
\\2A-\f p{2A}\mod {p^2}&\t{if $p=A^2+3B^2\e 1\mod 3$ and $3\mid A-1$.}
\endcases\endaligned$$
\endpro
Proof. Set $a=-\f 16$ and $b=-\f 12$. Then
$$\bp=\f{p-1}2\qtq{and}\ap=\cases \f{p-1}6&\t{if $p\e 1\mod 3$,}
\\\f{5p-1}6&\t{if $p\e 2\mod 3$.}
\endcases$$
For $p\e 2\mod 3$ we have $\ap>\bp$. Thus, by Theorem 4.1 we get
$$\aligned&\sum_{k=0}^{p-1}\b{-\f 16}k\b{-\f 13}k
\\&\e\f{(-1)^{\f{5p-1}6-\f{p-1}2-1}}{(\f{5p-1}6-\f{p-1}2)\b{\f{5p-1}6}{\f{p-1}2}}
\Big(-\f
12-\f{p-1}2\Big)=\f{3p}{2(p+1)\b{\f{5p-1}6}{\f{p-1}2}}\mod{p^2}.
\endaligned$$
Note that
$$\align \b{\f{5p-1}6}{\f{p-1}2}&=\b{\f{5p-1}6}{\f{p+1}3}
=\b{p-\f{p+1}6}{\f{p+1}3}\e \b{-\f{p+1}6}{\f{p+1}3}
\\&=\b{\f{p+1}6+\f{p+1}3-1}{\f{p+1}3}
=\b{\f{p-1}2}{\f{p+1}3}=\b{\f{p-1}2}{\f{p-5}6} \mod p.\endalign$$ We
then get the result in the case $p\e 2\mod 3$.
\par Now we assume $p\e 1\mod 3$ and so $p=A^2+3B^2\e 1\mod 3$ with $3\mid A-1$.
By Theorem 4.1 and Lemma 3.1 we obtain
$$\align&\sum_{k=0}^{p-1}\b{-\f 16}k\b{-\f 13}k
\\&\e \b{\f{p-1}2}{\f{p-1}6}\Big(1-\f p2H_{\f{p-1}2}+\f
p6H_{\f{p-1}6}+\f p3H_{\f{p-1}3}\Big)
\\&\e\b{\f{p-1}2}{\f{p-1}6}\Big(1-\f p2(-2q_p(2))+\f
p6\big(-2q_p(2)-\f 32q_p(3)\big)+\f p3\big(-\f 32q_p(3)\big)\Big)
\\&=\b{\f{p-1}2}{\f{p-1}6}\Big(1+p\big(\f 23q_p(2)-\f
34q_p(3)\big)\Big)\mod{p^2}.\endalign$$
 By [BEW,
Theorem 9.4.4],
$$\b{\f{p-1}2}{\f{p-1}6}
\e\Big(2A-\f p{2A}\Big)\Big(1-p\big(\f 23q_p(2)-\f
34q_p(3)\big)\Big)\mod{p^2}.$$ Therefore,
$$\align&\sum_{k=0}^{p-1}\b{-\f 16}k\b{-\f 13}k
\\&\e \Big(2A-\f p{2A}\Big)\Big(1^2-p^2\big(\f 23q_p(2)-\f
34q_p(3)\big)^2\Big) \e 2A-\f p{2A}\mod{p^2}.\endalign$$ This proves
the case $p\e 1\mod 3$. Hence the theorem is proved.

\pro{Conjecture 4.1} Let $p$ be an odd prime. Then
$$\sum_{k=0}^{p-1}\f{\b{-\f 14}k\b{-\f 12}k}{4^k}
\e\cases 2x-\f p{2x}\mod{p^2}&\t{if $p=x^2+y^2\e 1\mod{12}$ and
$2\nmid x$,}
\\2y-\f p{2y}\mod{p^2}&\t{if $p=x^2+y^2\e 5\mod{12}$ and
$2\nmid x$,}
\\0\mod{p^2}&\t{if $p\e 3\mod 4$.}\endcases$$\endpro

\pro{Conjecture 4.2} Let $p$ be an odd prime. Then
$$\aligned&\sum_{k=0}^{p-1}\f{\b{-\f 14}k\b{-\f 12}k}{(-3)^k}
\e (-1)^{\f{p-1}4}\sum_{k=0}^{p-1}\f{\b{-\f 14}k\b{-\f 12}k}{81^k}
\\&\e\cases 2x-\f p{2x}\mod {p^2}&\t{if $p=x^2+y^2\e 1\mod
4$ and $2\nmid x$,}
\\0\mod p&\t{if $p\e 3\mod 4$.}\endcases\endaligned$$\endpro

\pro{Conjecture 4.3} Let $p$ be an odd prime. Then
$$\aligned&\sum_{k=0}^{p-1}\f{\b{-\f 14}k\b{-\f 12}k}{(-80)^k}
\\&\e\cases 2x-\f p{2x}\mod {p^2}&\t{if $p=x^2+y^2\e \pm 1\mod
5$ and $2\nmid x$,}
\\2y-\f p{2y}\mod {p^2}&\t{if $p=x^2+y^2\e \pm 2\mod
5$ and $2\nmid x$,}
\\0\mod p&\t{if $p\e 3\mod 4$.}\endcases\endaligned$$\endpro

\pro{Conjecture 4.4} Let $p>5$ be a prime. Then
$$\aligned&\sum_{k=0}^{p-1}\b{-\f 14}k\b{-\f 12}k2^k
\\&=\cases 2x-\f p{2x}\mod {p^2}&\t{if $p=x^2+2y^2$ with $x\e 1\mod 4$,}
\\0\mod p&\t{if $ p\e 5,7\mod 8$.}
\endcases\endaligned$$
\endpro

\pro{Conjecture 4.5} Let $p>3$ be a prime. Then
$$\aligned&\sum_{k=0}^{p-1}\f{\b{-\f 16}k\b{-\f 13}k}{2^k}
\\&\e\cases \sls x3(2x-\f p{2x})\mod {p^2}&\t{if $p=x^2+6y^2\e 1,7\mod
{24}$,}\\\sls x3(2x-\f p{4x})\mod {p^2}&\t{if $p=2x^2+3y^2\e
5,11\mod {24}$,}
\\0\mod p&\t{if $p\e 13,17,19,23\mod {24}$.}\endcases\endaligned$$\endpro

\pro{Conjecture 4.6} Let $p>3$ be a prime. Then
$$\aligned\sum_{k=0}^{p-1}\f{\b{-\f 16}k\b{-\f 13}k}{(-16)^k}
\e\cases -\sls x3(x-\f px)\mod {p^2}&\t{if $4p=x^2+51y^2$,}
\\0\mod p&\t{if $\sls p{51}=-1$.}\endcases\endaligned$$\endpro

\pro{Conjecture 4.7} Let $p$ be a prime of the form $6k+5$. Then
$$\aligned\sum_{k=0}^{p-1}\f{\b{-\f 16}k\b{-\f 13}k}{(-80)^k}
\e \sum_{k=0}^{p-1}\f{\b{-\f 16}k\b{-\f 13}k}{(-3024)^k}\e 0\mod
p.\endaligned$$\endpro

\pro{Conjecture 4.8} Let $p>3$ be a prime. Then
$$\aligned\sum_{k=0}^{p-1}\f{\b{-\f 16}k\b{-\f 13}k}{(-1024)^k}
\e\cases -\sls x3(x-\f px)\mod {p^2}&\t{if $4p=x^2+123y^2$,}
\\0\mod p&\t{if $\sls p{123}=-1$.}\endcases\endaligned$$\endpro

\pro{Conjecture 4.9} Let $p>3$ be a prime. Then
$$\aligned\sum_{k=0}^{p-1}\f{\b{-\f 16}k\b{-\f 13}k}{(-250000)^k}
\e\cases -\sls x3(x-\f px)\mod {p^2}&\t{if $4p=x^2+267y^2$,}
\\0\mod p&\t{if $\sls p{267}=-1$.}\endcases\endaligned$$\endpro

\pro{Conjecture 4.10} Let $p>5$ be a prime. Then
$$\aligned\sum_{k=0}^{p-1}\b{-\f 13}k\b{-\f 16}k(-4)^k
=\cases \sls x3(2x-\f p{2x})\mod {p^2}&\t{if $p=x^2+15y^2$,}
\\ -\sls x3(10x-\f p{2x})\mod {p^2}&\t{if $p=5x^2+3y^2$,}
\\0\mod p&\t{if $\sls p{15}=-1$.}
\endcases\endaligned$$
\endpro

\pro{Conjecture 4.11} Let $p>3$ be a prime such that $p\e
13,17,19,23\mod {24}$. Then
$$\sum_{k=0}^{p-1}(-1)^k\b{-\f 16}k\b{-\f 23}k\e 0\mod p.$$
\endpro

\pro{Conjecture 4.12} Let $p$ be a prime such that $p\e 5\mod 6$.
Then
$$\sum_{k=0}^{p-1}\f{\b{-\f 16}k\b{-\f 23}k}{9^{2k}}
\e \sum_{k=0}^{p-1}\f{\b{-\f 16}k\b{-\f 23}k}{55^{2k}}\e 0\mod p.$$
\endpro

\pro{Conjecture 4.13} Let $p>5$ be a prime. Then
$$\aligned&\sum_{k=0}^{p-1}\b{-\f 12}k\b{-\f 13}k(-3)^k
\e\sum_{k=0}^{p-1}\f{\b{-\f 12}k\b{-\f 13}k}{(-27)^k} \\&\e \Ls
p5\sum_{k=0}^{p-1}\f{\b{-\f 12}k\b{-\f 13}k}{5^k}   \e
\Ls{-1}p\sum_{k=0}^{p-1}\b{-\f 12}k\b{-\f 13}k2^k
\\&\e\cases 2A-\f p{2A}\mod{p^2}&\t{if
$p=A^2+3B^2\e 1\mod 3$ with $3\mid A-1$,}\\0\mod p&\t{if $p\e 2\mod
3$.}
\endcases\endaligned$$\endpro

\pro{Conjecture 4.14} Let $p>5$ be a prime. Then
$$\aligned&\sum_{k=0}^{p-1}\f{\b{-\f 12}k\b{-\f 13}k}{(-4)^k}
\\&=\cases \sls p5 (2A-\f p{2A})\mod{p^2}\\\q\t{if
$p=A^2+3B^2\e 1\mod 3$ with $5\mid AB$ and $3\mid A-1$,} \\ \sls p5
(A+3B-\f p{A+3B})\mod{p^2}\\\q\t{if $p=A^2+3B^2\e 1\mod 3$ with
$A/B\e -1,-2\mod 5$ and $3\mid A-1$,}
\\0\mod
p\qq\t{if $p\e 2\mod 3$.}
\endcases\endaligned$$\endpro

\pro{Conjecture 4.15} Let $p\not=2,5$ be a prime such that
$\sls{-5}p=-1$. Then
$$\sum_{k=0}^{p-1}\b{-\f 18}k\b{-\f 58}k\f 1{5^k}
\e \sum_{k=0}^{p-1}\b{-\f 38}k\b{-\f 78}k\f 1{5^k}\e 0\mod p.$$
\endpro

\pro{Conjecture 4.16} Let $p\not=7$ be an odd prime such that
$\sls{-1}p=-1$. Then
$$\sum_{k=0}^{p-1}\b{-\f 18}k\b{-\f 58}k\f 1{49^k}
\e \sum_{k=0}^{p-1}\b{-\f 38}k\b{-\f 78}k\f 1{49^k}\e 0\mod p.$$
\endpro

\pro{Conjecture 4.17} Let $p\not=2,3$ be a prime such that
$\sls{-6}p=-1$. Then
$$\sum_{k=0}^{p-1}\b{-\f 18}k\b{-\f 58}k\f 1{(-8)^k}
\e \sum_{k=0}^{p-1}\b{-\f 38}k\b{-\f 78}k\f 1{(-8)^k}\e 0\mod p.$$
\endpro

\pro{Conjecture 4.18} Let $p>5$ be a prime such that $\sls{-2}p=-1$.
Then
$$\sum_{k=0}^{p-1}\b{-\f 18}k\b{-\f 58}k\f 1{(-2400)^k}\e
\sum_{k=0}^{p-1}\b{-\f 38}k\b{-\f 78}k\f 1{(-2400)^k}\e
 0\mod p.$$
\endpro

\pro{Conjecture 4.19} Let $p\not=2,5$ be a prime such that
$\sls{-10}p=-1$. Then
$$\sum_{k=0}^{p-1}\b{-\f 18}k\b{-\f 58}k\f 1{(-80)^k}\e
\sum_{k=0}^{p-1}\b{-\f 38}k\b{-\f 78}k\f 1{(-80)^k}\e 0\mod p.$$
\endpro

\pro{Conjecture 4.20} Let $p\not=2,3,19$ be a prime with $\sls
{-19}p=-1$. Then
$$\sum_{k=0}^{p-1}\b{-\f 1{12}}k\b{-\f 7{12}}k\f 1{513^k}\e 0\mod p.$$
\endpro

\pro{Conjecture 4.21} Let $p\not=2,3,17$ be a prime with $\sls
{-51}p=-1$. Then
$$\sum_{k=0}^{p-1}\b{-\f 1{3}}k\b{-\f 56}k\f 1{17^k}\e 0\mod p.$$
\endpro

\pro{Conjecture 4.22} Let $p\not=2,3$ be a prime with $\sls
{-3}p=-1$. Then
$$\sum_{k=0}^{p-1}\b{-\f 1{3}}k\b{-\f 56}k\f 1{81^k}\e 0\mod p.$$
\endpro

\pro{Conjecture 4.23} Let $p\not=2,3$ be a prime with $\sls
{-6}p=-1$. Then
$$\sum_{k=0}^{p-1}\b{-\f 1{3}}k\b{-\f 56}k(-1)^k\e 0\mod p.$$
\endpro

\subheading{5. Congruences for $\sum_{k=0}^{p-1}\b{-\f 13}k^2m^k$
and $\sum_{k=0}^{p-1}\b{-\f 14}k^2m^k$}

\pro{Theorem 5.1} Let $p$ be an odd prime. Then
$$\aligned&\sum_{k=0}^{p-1}\b{-\f 14}k^2
\\&\e \cases 2p\b{\f{(p-1)}2}{\f{p+1}4}^{-1}\mod{p^2}&\t{if $p\e
3\mod 4$,}
\\2x-\f p{2x}\mod{p^2}&\t{if $p=x^2+y^2\e 1\mod 4$
and $4\mid x-1$.}
\endcases\endaligned$$
\endpro
Proof. Set $a=-\f 14$ and $b=-\f 12$. Then $\bp=\f{p-1}2$, $b-a=-\f
14$ and
$$\ap=\langle b-a\rangle_p=\cases \f{p-1}4&\t{if $p\e 1\mod 4$,}
\\\f{3p-1}4&\t{if $p\e 3\mod 4$.}\endcases$$

If $p\e 3\mod 4$, then $\ap>\bp$. Thus, by Theorem 4.1 we have
$$\aligned\sum_{k=0}^{p-1}\b{-\f 14}k^2
\e \f{(-1)^{\f{3p-1}4-\f{p-1}2-1}}
{(\f{3p-1}4-\f{p-1}2)\b{\f{3p-1}4}{\f{p-1}2}} \Big(-\f p2 \Big) \e
2p(-1)^{\f{p+1}4}\b{\f{3p-1}4}{\f{p-1}2}^{-1}\mod{p^2}.\endaligned$$
Since
$$\align \b{\f{3p-1}4}{\f{p-1}2}&=\b{\f{3p-1}4}{\f{p+1}4}
=\b{p-\f{p+1}4}{\f{p+1}4} \e
\b{-\f{p+1}4}{\f{p+1}4}\\&=(-1)^{\f{p+1}4}\b{\f{p+1}4
+\f{p+1}4-1}{\f{p+1}4} =(-1)^{\f{p+1}4}\b{\f{p-1}2}{\f{p+1}4}\mod
p,\endalign$$ by the above we obtain the result in the case $p\e
3\mod 4$.

\par Now we assume $p\e 1\mod 4$ and so $p=x^2+y^2$ with $x,y\in\Bbb
Z$ and $x\e 1\mod 4$. By Lemma 3.1(i) we have $H_{\f{p-1}4}\e
-3q_p(2)\mod p$ and $H_{\f{p-1}2}\e -2q_p(2)$. Now applying the
above and Theorem 4.1 we deduce
$$\aligned\sum_{k=0}^{p-1}\b{-\f 14}k^2
&\e\b{\f{p-1}2}{\f{p-1}4}\Big(1-\f p2 H_{\f{p-1}2}+ \f
p2H_{\f{p-1}4}\Big)
\\&\e \b{\f{p-1}2}{\f{p-1}4}\Big(1-\f p2(-2q_p(2))+\f
p2(-3q_p(2))\Big)
\\&=\b{\f{p-1}2}{\f{p-1}4}\Big(1-\f 12pq_p(2)\Big)\mod{p^2}.
\endaligned$$
By [BEW] we have
$$\b{\f{p-1}2}{\f{p-1}4}\e \Big(2x-\f
p{2x}\Big)\Big(1+\f 12pq_p(2)\Big)\mod{p^2}.$$ Hence
$$\sum_{k=0}^{p-1}\b{-\f 14}k^2\e 2x-\f
p{2x}\mod{p^2}.$$ This proves the result in the case $p\e 1\mod 4$.
The proof is now complete.
\par\q
\newline{\bf Remark 5.1} Let $p$ be an odd prime and $a\in\Bbb Z_p$.
From [Su4] we know that
$$\sum_{k=0}^{p-1}\b ak^2\e \b{2a}{\ap}\mod{p^2}.$$
Hence
$$\sum_{k=0}^{p-1}\b{-\f 14}k^2\e \cases \b{-\f
12}{\f{p-1}4}\mod{p^2}&\t{if $p\e 1\mod 4$,}
\\\b{-\f 12}{\f{3p-1}4}\mod{p^2}&\t{if $p\e 3\mod 4$.}
\endcases$$

\pro{Lemma 5.1 ([G, (3.134)])} Let $n$ be a nonnegative integer and
$x\not=1$. Then
$$P_n(x)=\Ls{x-1}2^n\sum_{k=0}^n\b nk^2\Ls{x+1}{x-1}^k.$$
\endpro
\pro{Lemma 5.2} Let $p$ be an odd prime and $a\in\Bbb Z_p$ with
$a\not\e 0\mod p$. Then
$$\sum_{k=0}^{p-1}\b{a}k^2t^k\e t^{\ap}
\sum_{k=0}^{p-1}\f{\b{a}k^2}{t^k}\mod p$$ and
$$\sum_{k=0}^{p-1}\b{a}k^2t^k\e (t-1)^{\ap}P_{\ap}
\Ls{t+1}{t-1}\mod p.$$
\endpro
Proof. It is clear that
$$\align\sum_{k=0}^{p-1}\b{a}k^2t^k
&\e \sum_{k=0}^{\ap}\b{\ap}k^2t^k
=\sum_{s=0}^{\ap}\b{\ap}{\ap-s}^2t^{\ap-s}
\\&=t^{\ap}\sum_{k=0}^{\ap}\b{\ap}k^2t^{-k}
\e t^{\ap} \sum_{k=0}^{p-1}\f{\b{a}k^2}{t^k}\mod p.\endalign$$ Using
Lemma 5.1 we see that for $x\not=1$,
$$\align P_{\ap}(x)&=\Ls{x-1}2^{\ap}\sum_{k=0}^{\ap}
\b {\ap}k^2\Ls{x+1}{x-1}^k\\& =\Ls{x-1}2^{\ap}\sum_{k=0}^{p-1} \b
{\ap}k^2\Ls{x+1}{x-1}^k
\\&\e \Ls{x-1}2^{\ap}\sum_{k=0}^{p-1} \b
{a}k^2\Ls{x+1}{x-1}^k\mod p.\endalign$$ Set $x=\f{t+1}{t-1}$. Then
$t=\f{x+1}{x-1}$ and $\f{x-1}2=\f 1{t-1}$. Now substituting $x$ with
$\f{t+1}{t-1}$ in the above congruence we obtain the remaining
result. \pro{Lemma 5.3} Let $p$ be an odd prime and
$m\in\{1,2,\ldots,\f{p-1}2\}$. Then $P_{p-1-m}(x)\e P_m(x)$.
\endpro
Proof. Since $m<\f p2$ we have $p-1-m\ge m$. It is well-known that
$\b{-t}k=(-1)^k\b{t+k-1}k$ and
$$P_n(x)=\sum_{k=0}^n\b nk\b{n+k}k\Ls{x-1}2^k.$$
Thus,
$$\align P_{p-1-m}(x)&=\sum_{k=0}^{p-1-m}
\b{p-1-m}k\b{p-1-m+k}k\Ls{x-1}2^k \\&\e \sum_{k=0}^{p-1-m}
\b{-1-m}k\b{-1-m+k}k\Ls{x-1}2^k
\\&=\sum_{k=0}^{p-1-m}(-1)^k\b{m+k}k\cdot (-1)^k\b mk\Ls{x-1}2^k
\\&=\sum_{k=0}^m\b mk\b{m+k}k\Ls{x-1}2^k
=P_m(x)\mod p.\endalign$$ This proves the lemma.

\pro{Theorem 5.2} Let $p$ be an odd prime and $a,t\in\Bbb Z_p$ with
$a\not\e -1\mod p$ and $t\not\e 1\mod p$. Then
$$\sum_{k=0}^{p-1}\b {-1-a}k^2t^k\e (t-1)^{-2\ap}
\sum_{k=0}^{p-1}\b{a}k^2t^k\mod p.$$
\endpro
Proof. Clearly $\ag{-1-a}=p-1-\ap$. Thus, using Lemmas 5.2 and 5.3
we see that for $t\not\e 1\mod p$,
$$\align &\sum_{k=0}^{p-1}\b {-1-a}k^2t^k
\\&\e (t-1)^{\ag{-1-a}}P_{\ag{-1-a}}
\Ls{t+1}{t-1}=(t-1)^{p-1-\ap}P_{p-1-\ap}\Ls{t+1}{t-1}
\\&\e (t-1)^{-\ap}P_{\ap}\Ls{t+1}{t-1}
\e (t-1)^{-2\ap}\sum_{k=0}^{p-1}\b ak^2t^k \mod p.\endalign$$ Tis
proves the theorem.

\pro{Theorem 5.3} Let $p>3$ be a prime. Then
$$\aligned&\sum_{k=0}^{p-1}\b{-\f 13}k^29^k\e \f 1{3^{[\f p3]}}
\sum_{k=0}^{p-1}\b{-\f 13}k^2\f 1{9^k}
\\&\e\cases L\mod p&\t{if $p\e 1\mod 3$, $4p=L^2+27M^2$ and
$3\mid L-2$,}
\\0\mod{p}&\t{if $p\e 2\mod 3$.}\endcases\endaligned$$\endpro
Proof. By Lemmas 5.2 and 5.3 we have
$$\aligned &9^{\langle -\f 13\rangle_p}
\sum_{k=0}^{p-1}\b{-\f 13}k^2\f 1{9^k}\\&\e \sum_{k=0}^{p-1}\b{-\f
13}k^29^k \e 8^{\langle -\f 13\rangle_p}P_{\langle -\f
13\rangle_p}\Ls{10}8
\\&=\cases 8^{\f{p-1}3}P_{\f{p-1}3}\sls 54 \e P_{\f{p-1}3}\sls 54\mod
p &\t{if $p\e 1\mod 3$,}
\\8^{\f{2p-1}3}P_{\f{2p-1}3}\sls 54\e 2P_{\f{p-2}3}\sls 54\mod
p&\t{if $p\e 2\mod 3$.}\endcases\endaligned$$ From [S4, Theorem 3.2]
we know that
$$P_{[\f p3]}\Ls 54\e \cases L\mod p&\t{if $p\e 1\mod 3$,}
\\0\mod p&\t{if $p\e 2\mod 3$.}\endcases$$
Since
$$9^{\ag{-\f 13}}=\cases 9^{\f{p-1}3}\e 3^{-\f{p-1}3}\mod p
&\t{if $p\e 1\mod 3$,}
\\9^{\f{2p-1}3}\e 3^{-\f{p-2}3}\mod p&\t{if $p\e 2\mod 3$,}
\endcases$$ combining all the above we deduce the result.

 \pro{Conjecture 5.1} Let $p$ be a prime such that $p\e 1\mod 3$,
  $4p=L^2+27M^2(L,M\in\Bbb Z)$ and $L\e 2\mod 3$.  Then
$$\sum_{k=0}^{p-1}\b{-\f 13}k^29^k
\e L-\f pL\mod{p^2}$$\endpro

\pro{Theorem 5.4} Let $p$ be an odd prime. Then
$$\aligned&\sum_{k=0}^{p-1}\b{-\f 14}k^2(-8)^k
\\&\e \cases (-1)^{\f{p-1}4}2x\mod p&\t{if $p=x^2+y^2\e 1\mod 4$
and $4\mid x-1$,}
\\0\mod p&\t{if $p\e 3\mod 4$}\endcases\endaligned$$
and
$$\aligned&\sum_{k=0}^{p-1}\b{-\f 14}k^2\f 1{(-8)^k}
\\&\e\cases (-1)^{\f y4}2x\mod p
&\t{if $p=x^2+y^2\e 1\mod 8$ and $4\mid x-1$,}\\(-1)^{\f{y-2}4}
2y\mod p&\t{if $p=x^2+y^2\e 5\mod 8$ and $2\nmid x$,}
\\0\mod p&\t{if $p\e 3\mod 4$.}\endcases\endaligned$$
\endpro
Proof. By Lemmas 5.2 and 5.3 we have
$$\aligned&\sum_{k=0}^{p-1}\b{-\f 14}k^2(-8)^k
\\&\e (-8)^{\ag{-\f 14}}\sum_{k=0}^{p-1}\b{-\f 14}k^2\f 1{(-8)^k}
\\&\e (-9)^{\ag{-\f 14}}P_{\ag{-\f 14}}\Ls 79=9^{\ag{-\f 14}}P_{\ag{-\f
14}}\Big(-\f 79\Big)
\\&=\cases 9^{\f{p-1}4}P_{\f{p-1}4}(-\f 79)\e \sls 3p
P_{\f{p-1}4}(-\f 79)\mod p&\t{if $p\e 1\mod 4$,}
\\9^{\f{3p-1}4}P_{\f{3p-1}4}(-\f 79)\e 9^{\f{3p-1}4}
P_{\f{p-3}4}(-\f 79)\mod p&\t{if $p\e 3\mod
4$.}\endcases\endaligned$$ From [S4, Theorem 2.4] we know that
$$P_{[\f p4]}\Big(-\f 79\Big)\e \cases (-1)^{\f{p-1}4}\sls p32x\mod
p&\t{if $p=x^2+y^2\e 1\mod 4$ and $4\mid x-1$,}
\\0\mod p&\t{if $p\e 3\mod 4$.}
\endcases$$
When $p\e 1\mod 4$, we have   $(-8)^{-\ag{-\f
14}}=(-8)^{-\f{p-1}4}\e (-2)^{\f{p-1}4}\mod p.$ It is well known
that (see [BEW])
$$2^{\f{p-1}4}\e \cases(-1)^{\f y4}\mod p&\t{if $p\e 1\mod 8$,}
\\(-1)^{\f{y-2}4}\f yx\mod p&\t{if $p\e 5\mod 8$.}
\endcases$$
Now combining all the above we deduce the result.

\pro{Conjecture 5.2} Let $p$ be an odd prime. Then
$$\aligned&\sum_{k=0}^{p-1}\b{-\f 14}k^2(-8)^k
\\&\e\cases (-1)^{\f{p-1}4}(2x-\f p{2x})\mod {p^2}&\t{if $p=x^2+y^2\e 1\mod
4$ and $4\mid x-1$,}
\\0\mod p&\t{if $p\e 3\mod 4$.}\endcases\endaligned$$
and
$$\aligned&\sum_{k=0}^{p-1}\f{\b{-\f 14}k^2}{(-8)^k}
\\&\e\cases (-1)^{\f y4}(2x-\f p{2x})\mod {p^2}
&\t{if $p=x^2+y^2\e 1\mod 8$ and $4\mid x-1$,}\\(-1)^{\f{y-2}4}
(2y-\f p{2y})\mod {p^2}&\t{if $p=x^2+y^2\e 5\mod 8$ and $4\mid
y-2$,}
\\0\mod p&\t{if $p\e 3\mod 4$.}\endcases\endaligned$$\endpro

\pro{Theorem 5.5} Let $p>3$ be a prime. Then
$$\aligned&\sum_{k=0}^{p-1}\b{-\f 14}k^24^k
\e \f{3-(-1)^{\f{p-1}2}}2\Ls 2p\sum_{k=0}^{p-1}\b{-\f 14}k^2\f
1{4^k}
\\&\e \cases (-1)^{\f{p-1}4+\f{A-1}2}2A\mod p&\t{if $p=A^2+3B^2
\e 1\mod {12}$,}
\\(-1)^{\f{p+1}4}6B\mod p&\t{if $p=A^2+3B^2\e 7\mod {12}$
and $4\mid B-1$,}
\\0\mod p&\t{if $p\e 2\mod 3$}\endcases\endaligned$$
\endpro
Proof. By Lemmas 5.2 and 5.3 we have
$$\aligned&\sum_{k=0}^{p-1}\b{-\f 14}k^24^k
\\&\e 4^{\ag{-\f 14}}\sum_{k=0}^{p-1}\b{-\f 14}k^2\f 1{4^k}
\\&\e 3^{\ag{-\f 14}}P_{\ag{-\f 14}}\Ls 53=(-3)^{\ag{-\f 14}}P_{\ag{-\f
14}}\Big(-\f 53\Big)
\\&=\cases (-3)^{\f{p-1}4}P_{\f{p-1}4}(-\f 53)\mod p&\t{if $p\e 1\mod 4$,}
\\(-3)^{\f{3p-1}4}P_{\f{3p-1}4}(-\f 53)\e (-3)^{-\f{p-3}4}
P_{\f{p-3}4}(-\f 53)\mod p&\t{if $p\e 3\mod
4$.}\endcases\endaligned$$ From [S4, Theorem 2.5] we know that
$$P_{[\f p4]}\Big(-\f 53\Big)\e \cases 2A\mod
p&\t{if $p=A^2+3B^2\e 1\mod 3$ and $3\mid A-1$,}
\\0\mod p&\t{if $p\e 2\mod 3$.}
\endcases$$
Hence the result is true for $p\e 2\mod 3$. \par Now assume
$p=A^2+3B^2\e 1\mod 3$ and $A\e 1\mod 3$. If $p\e 1\mod {12}$, by
[S5, p.1317] we have $3^{\f{p-1}4}\e (-1)^{\f{A-1}2}\mod p$ and
$4^{\ag{-\f 14}}=4^{\f{p-1}4}\e \ls 2p=(-1)^{\f{p-1}4}\mod p$. Hence
$$\sum_{k=0}^{p-1}\b{-\f 14}k^24^k
\e (-1)^{\f{p-1}4}\sum_{k=0}^{p-1}\b{-\f 14}k^2\f 1{4^k} \e
(-1)^{\f{p-1}4+\f{A-1}2}2A\mod p.$$ If $p\e 7\mod {12}$ and $B\e
1\mod 4$, by [S5, p.1317] we have $3^{\f{p-3}4}\e \f BA\mod p$.
Since $4^{\ag{-\f 14}}=4^{\f{3p-1}4}=2^{p-1+\f{p+1}2} \e 2\sls
2p\mod p$, by the above we get
$$\sum_{k=0}^{p-1}\b{-\f 14}k^24^k
\e 2\Ls 2p\sum_{k=0}^{p-1}\b{-\f 14}k^2\f 1{4^k} \e
(-1)^{\f{p-3}4}\f AB\cdot 2A\e (-1)^{\f{p+1}4}6B\mod p.$$ Now
combining all the above we deduce the result.

\pro{Conjecture 5.3} Let $p>3$ be an odd prime. Then
$$\aligned&\sum_{k=0}^{p-1}\b{-\f 14}k^24^k\e
\sum_{k=0}^{p-1}\b{-\f 14}k\b{-\f 12}k(-8)^k
\\&\e\cases (-1)^{\f{p-1}4+\f{A-1}2}
(2A-\f p{2A})\mod {p^2}&\t{if $p=A^2+3B^2\e 1\mod {12}$,}
\\(-1)^{\f{p+1}4+\f{B-1}2}
(6B-\f p{2B})\mod {p^2}&\t{if $p=A^2+3B^2\e 7\mod {12}$,}
\\0\mod p&\t{if $p\e 2\mod 3$}\endcases\endaligned$$
and
$$\aligned\sum_{k=0}^{p-1}\f{\b{-\f 14}k^2}{4^k}
\e\cases (-1)^{\f{A-1}2}(2A-\f p{2A})\mod {p^2}&\t{if $p=A^2+3B^2\e
1\mod{12}$,}
\\(-1)^{\f{B-1}2}(3B-\f p{4B})\mod {p^2}&\t{if $p=A^2+3B^2\e
7\mod{12}$,}
\\0\mod p&\t{if $p\e 2\mod 3$.}\endcases\endaligned$$\endpro

\pro{Theorem 5.6} Let $p\not=2,7$ be a prime. Then
$$\aligned\sum_{k=0}^{p-1}\b{-\f 14}k^264^k
&\e\f{9-7(-1)^{\f{p-1}2}}2\Ls 2p\sum_{k=0}^{p-1}\b{-\f 14}k^2\f
1{64^k}
\\&\e\cases
(-1)^{\f{p-1}4+\f{x-1}2}2x\mod p&\t{if $p=x^2+7y^2\e 1\mod
4$,}\\(-1)^{\f{p+1}4+\f{y-1}2}42y\mod p&\t{if $p=x^2+7y^2\e 3\mod
4$,}
\\0\mod p&\t{if $p\e 3,5,6\mod 7$.}\endcases\endaligned$$\endpro

Proof. By Lemmas 5.2 and 5.3 we have
$$\aligned&\sum_{k=0}^{p-1}\b{-\f 14}k^264^k
\\&\e 64^{\ag{-\f 14}}\sum_{k=0}^{p-1}\b{-\f 14}k^2\f 1{64^k}
\\&\e 63^{\ag{-\f 14}}P_{\ag{-\f 14}}\Ls {65}{63}
=(-63)^{\ag{-\f 14}}P_{\ag{-\f 14}}\Big(-\f {65}{63}\Big)
\\&=\cases (-63)^{\f{p-1}4}P_{\f{p-1}4}(-\f {65}{63})\mod p&\t{if $p\e 1\mod 4$,}
\\(-63)^{\f{3p-1}4}P_{\f{3p-1}4}(-\f {65}{63})\e (-63)^{-\f{p-3}4}
P_{\f{p-3}4}(-\f {65}{63})\mod p&\t{if $p\e 3\mod
4$.}\endcases\endaligned$$ By [S4, Theorem 2.6],
$$P_{[\f p4]}\Big(-\f{65}{63}\Big)
\e \cases 2x\sls p3\sls x7\mod p&\t{if $p=x^2+7y^2\e 1,2,4\mod 7$,}
\\0\mod p&\t{if $p\e 3,5,6\mod 7$.}
\endcases$$
Hence the result is true for $p\e 3,5,6\mod 7$.
\par Now suppose $p\e 1,2,4\mod 7$ and so $p=x^2+7y^2$ with $x,y\in
\Bbb Z$. If $p\e 1\mod 4$, by [S5, p.1317] we have $7^{\f{p-1}4} \e
(-1)^{\f{x-1}2}\sls x7\mod p$ and so
$$\align&(-63)^{\f{p-1}4}P_{\f{p-1}4}(-\f {65}{63})
\\&\e (-1)^{\f{p-1}4}\Ls 3p
\cdot (-1)^{\f{x-1}2}\Ls x7\cdot 2x\Ls p3\Ls x7
=(-1)^{\f{p-1}4+\f{x-1}2} 2x\mod p.\endalign$$ If $p\e 3\mod 4$, by
[S5, p.1317] we have $7^{\f{p-3}4} \e (-1)^{\f{y+1}2}\sls x7\f
yx\mod p$ and so
$$\align&(-63)^{-\f{p-3}4}P_{\f{p-3}4}(-\f {65}{63})
\\&\e (-1)^{\f{p-3}4}3\Ls 3p
\cdot (-1)^{\f{y+1}2}\Ls x7\f xy\cdot 2x\Ls p3\Ls x7 \e
(-1)^{\f{p+1}4+\f{y-1}2} 42y\mod p.\endalign$$ Note that
$$64^{\ag{-\f 14}}=\cases 64^{\f{p-1}4}\e \sls 2p\mod p&\t{if $p\e
1\mod 4$,}
\\64^{\f{3p-1}4}\e 8\ls 2p\mod p&\t{if $p\e 3\mod 4$.}
\endcases$$
Combining all the above we deduce the result.

\pro{Conjecture 5.4} Let $p\not=2,7$ be a prime. Then
$$\aligned\sum_{k=0}^{p-1}\b{-\f 14}k^264^k
\e\cases \sls 2p(-1)^{\f{x-1}2}(2x-\f p{2x})\mod {p^2}&\t{if
$p=x^2+7y^2\e 1\mod 4$,}\\\sls 2p(-1)^{\f{y-1}2}(42y-\f
{3p}{2y})\mod {p^2}&\t{if $p=x^2+7y^2\e 3\mod 4$,}
\\0\mod p&\t{if $p\e 3,5,6\mod 7$.}\endcases\endaligned$$
and
$$\aligned\sum_{k=0}^{p-1}\f{\b{-\f 14}k^2}{64^k}
\e\cases (-1)^{\f{x-1}2}(2x-\f p{2x})\mod {p^2}&\t{if $p=x^2+7y^2\e
1\mod 4$,}\\\f{3}4(-1)^{\f{y-1}2}(7y-\f p{4y})\mod {p^2}&\t{if
$p=x^2+7y^2\e 3\mod 4$,}
\\0\mod p&\t{if $p\e 3,5,6\mod 7$.}\endcases\endaligned$$\endpro

\pro{Conjecture 5.5} Let $p$ be a prime such that $p\e 5,7\mod 8$.
Then
$$\sum_{k=0}^{p-1}(-1)^k\b{-\f 14}k^2\e
\cases (-1)^{\f{x+1}2}(2x-\f p{2x})\mod {p^2} &\t{if $p=x^2+2y^2\e
1\mod 8$,}\\(-1)^{\f{y-1}2}(4y-\f p{2y})\mod {p^2} &\t{if
$p=x^2+2y^2\e 3\mod 8$},\\0\mod p&\t{if $p\e 5,7\mod 8$.}\endcases$$
\endpro

\enddocument
\bye